\documentclass[12pt,centertags,tbtags]{amsart}
\usepackage{xspace}
\usepackage{geometry}
\newcommand\R{\ensuremath{\mathbb{R}}\xspace}  
 
\newcommand{\aop}{{\mathbf a}}  
\newcommand\indi[1]{\nbOne_{#1}}
\newcommand\nbOne{{\mathchoice {\mathrm {1\mskip-4.1mu l}} {\mathrm{ 1\mskip-4.1mu
        l}} {\mathrm {1\mskip-4.6mu l}} {\mathrm {1 \mskip-5.2mul}}}}

\newcommand\diw{\operatorname{div}}
\makeatletter
\def\th@definition{%
  \normalfont 
}
\def\th@plain{%
  \slshape 
}
\def\th@remark{%
  \normalfont 
  \thm@preskip\topsep
  \divide\thm@preskip\tw@
  \thm@postskip\thm@preskip
}
\makeatother
\theoremstyle{remark}
\newtheorem{remark}{Remark}[section]

\theoremstyle{definition}
\newtheorem{definition}[remark]{Definition}
\theoremstyle{plain}
\newtheorem{proposition}[remark]{Proposition}
\newtheorem{lemma}[remark]{Lemma}
\newtheorem{theorem}[remark]{Theorem}

\numberwithin{equation}{section}

\address{
Laboratoire de~Math\'ematiques~Rapha\"el~Salem
UMR 6085 CNRS -- Universit\'e de Rouen, 
Avenue de l'Universit\'e, BP.12, 
F76801 Saint-\'Etienne du Rouvray}
\email{ Olivier.Guibe@univ-rouen.fr}

\author[{M. Ben Cheikh Ali, O. Guib\smash{\'e}}]{M. Ben Cheikh Ali, O. Guib\'e}

\title
{  Nonlinear and non-coercive elliptic problems with integrable data
}

\begin{document}

\maketitle
\begin{abstract}
In this paper we study existence and uniqueness of renormalized solution
to the following problem%
\begin{equation*}
\left\{ 
\begin{array}{l}
\lambda \left( x,u\right) -\diw\left( \aop\left( x,Du\right) +%
\Phi \left( x,u\right) \right) =f\text{ \ \ \ in }\Omega ,
\\ 
\aop \left( x,Du\right) +\Phi \left( x,u\right) \cdot \mathbf{n%
}=0\text{ \ on }\Gamma _{n}, \\ 
u=0\text{ \ on }\Gamma _{d}.%
\end{array}%
\right. 
\end{equation*}
The main difficulty in this task is that in general the operator
entering in
the above equation is not coercive in a Sobolev space. Moreover, the
possible degenerate character of $\lambda$ with respect to $u$ renders
more complex the proof of uniqueness for integrable data $f$.

\end{abstract}

\section{Introduction}

In the present paper we study the class of nonlinear equations of the type 
\begin{eqnarray}
\lambda ( x,u) -\diw ( \aop ( x,Du) + 
\Phi ( x,u) )  &=& f \text{ \ \ \ in } 
\Omega ,  \label{1.1} \\
\big(\aop ( x,Du) +\Phi ( x,u)\big) \cdot \mathbf{n 
} &=&0\text{ \ on }\Gamma _{n},  \label{1.2} \\
u &=&0\text{ \ on }\Gamma _{d},  \label{1.3}
\end{eqnarray} 
where $\Omega $ is a bounded connected and open subset of $\mathbb{R}^{N}$ ($N\geq 2$) with
Lipschitz boundary $\partial \Omega $,   $\Gamma _{n}$ and $\Gamma _{d}$ are
such that $\Gamma _{n}\cup \Gamma _{d}=\partial \Omega $, $\Gamma _{n}\cap \Gamma
_{d}$ =$\emptyset $ and $\sigma ( \Gamma _{d}) >0$ (where $\sigma 
$ denotes the $N-1$ dimensional Lebesgue--measure on $\partial \Omega $). The
vector $\mathbf{n}$ is the outer unit normal to $\partial \Omega $ and the
data $ f $ is assumed to belong to $L^{1}( \Omega ) $. The
operator $u\mapsto -\diw ( \aop ( x,Du) ) $ is
monotone (but not necessarily strictly monotone) from the Sobolev space $ 
W_{0}^{1,p}( \Omega ) $ into $W^{-1,p^{\prime }}( \Omega
) $ with $1<p\leq N$ ($p^{\prime }=p/(p-1)$). The
functions  $\lambda\,:\,\Omega \times \mathbb{ 
R\longmapsto R} $ and $ 
\Phi\, :\,\Omega \times \mathbb{R\longmapsto R}^{N} $ are
Carath\'{e}odory functions such that
$\lambda ( x,r) r\geq 0$ for any $r\in\R$, almost
everywhere in $\Omega$ and such that
$|\Phi ( x,r) | \leq {b}( x)
( 1+|r|) ^{p-1}$ for any $r\in\R$, almost everywhere  in $\Omega $ with
$ b $ satisfying
some appropriate summability hypotheses that depend on $p$ and $N$ (see
condition (\ref{2.5bis}) below).
\par
Problem (\ref{1.1})--(\ref{1.3}) is motivated by the homogenization in the
particular case where $\aop  
( x,\xi)=A(x)\xi$  and where $\Omega $ is a perforated domain
with Neumann condition  on the boundary of the holes 
and Dirichlet  condition on the outside boundary of $ 
\Omega $ (see  \cite{Be99} and \cite{Moh2001}).
\par

 The main difficulty in dealing with
the existence of a solution of (\ref{1.1})--(\ref{1.3}) is the lack of
coercivity due to the term $-\diw(\Phi(x,u))$. 
As an example, consider the pure
Dirichlet case (i.e. $\Gamma_n= \emptyset$), the operator  $\aop  
( x,Du) +\Phi ( x,u) =\left| Du\right|
^{p-2}Du+ b ( x) \left| u\right| ^{p-2}u$ with $ b  
\in L^{{N}/(p-1)}( \Omega ) $. Then thanks to Sobolev's embedding
theorem, the  operator $u\longmapsto - 
\diw ( \aop ( x,Du) +\Phi (
x,u) ) $ is well defined from $W_{0}^{1,p}( \Omega ) $
to $W^{-1,p^{\prime }}( \Omega ) $ but it is not  coercive in
general except if $\left\|  b \right\| _{L^{N/(p-1)}(
\Omega ) }$ is small enough.

\par
Existence results for some similar non-coercive problems (with in addition 
lower order terms) are proved  in \cite{DelPos98} when $ f \in W^{-1,p^{\prime }}( \Omega
)$ and in 
 \cite{DelPos96} and \cite{Liang92} when $f$ is a Radon measure with
 bounded total variation in $\Omega$ (solutions in the sense of
 distributions are then used in this case).
A non-coercive linear case is studied  in \cite{droniou02}. 
In \cite{droniou03} the author gives local and global estimates
 for nonlinear non-coercive equations with measure data (with a
 stronger assumption of type (\ref{2.5}) below than the one used in the
 present paper, see (\ref{2.5bis}), in the case $p=N$).
  Entropy solutions to similar equations are considered in \cite{B}.

For integrable data $f$ we give in the present paper an existence
result (see Theorem \ref{th3.1} in Section 3) using the framework of renormalized solution. This
notion has been introduced by R. J. DiPerna and P.-L. Lions in \cite{DPL89a} for first
order equations and has been developed for elliptic problems with $L^{1}$
data in \cite{Mur93} (see also \cite{LM}). In
\cite{DMOP99} the authors give a definition of a renormalized solution
 for elliptic problems
with general measure data
and prove the existence of such a solution (a class of  nonlinear
elliptic equations with lower-order terms which are not coercive and
right-hand side measure is also studied in \cite{BMMP99}).

\par
\smallskip
 Another interesting question related to problem (\ref{1.1})--(\ref 
{1.3}) deals with the uniqueness of a solution. In \cite{Mur93} F. Murat
proves that the renormalized solution of 
\begin{eqnarray*}
\lambda u-\diw ( A ( x) Du+\phi  
( u) )  &=& f \text{ \ \ \ in }\Omega , \\
u &=&0\text{ on }\partial \Omega,
\end{eqnarray*} 
where $ f\in L^{1}( \Omega)$ and $\lambda >0$ is unique as soon as
$\phi $ is a locally Lipschitz continuous vector field.
In this result it is important to assume that  $\phi $ does not depend on 
$x$ together with pure homogeneous Dirichlet boundary conditions (see also \cite{LM}
and \cite{rak94} for more general operators and \cite{BR98} in the parabolic
case). When $\lambda ( x,s) $ is strictly monotone we prove in
Theorem \ref{th4.1} that the renormalized solution of (\ref{1.1})--(\ref{1.3}) is
unique if $\Phi(x,s) $ is locally Lipschitz continuous with respect to the
second variable. 
\par
As far as the case $\lambda ( x,u) \equiv 0$ is
concerned, gathering the result of \cite{Ar86}, \cite{BGM92}
 and \cite{CM89}, let us recall that when
$1<p\leq 2$ and  
$f\in W^{-1,p^{\prime }}( \Omega ) $, the uniqueness of
the variational solution of 
\begin{equation*}
\left\{ 
\begin{array}{l}
-\diw ( \aop ( x,u,Du) +\phi (
u) ) = f \text{ \ \ in }\Omega , \\ 
u=0\text{ on }\partial \Omega , 
\end{array} 
\right. 
\end{equation*} 
is obtained under strongly monotonicity assumption on the operator $\aop ( x,s,\xi)  $
and under  global Lipschitz
conditions on the functions $\aop  
( x,s,\xi ) $ and $\phi ( s) $ 
with respect to the variable $s$ (or a strong
control of the modulus of continuity). Moreover uniqueness may fail
 if  $2<p<\infty $ (see \cite{BGM92}). In the quasi-linear case
 (i.e. $a(x,s,\xi)=A(x,s)\xi$) and for integrable data uniqueness
 results have been  obtained in  \cite 
{Po01}  under a very general condition on the
 matrix field $A$ and the function $\phi$ (the author uses strongly the quasi-linear character
 of the problem). 
When $\lambda(x,s)\equiv 0$  we investigate in the present paper the uniqueness question
in the nonlinear case 
for  $1<p\leq 2$ and integrable data; global conditions on $\aop$ and
 $\Phi$ which insure uniqueness of the renormalized solution are given in Theorem \ref{th4.2}.

\par

The content of the paper is as follows.  In section 2 we precise the 
assumptions on the data and we give the definition of a renormalized solution of
problem (\ref{1.1})--(\ref{1.3}). This section is completed  by giving
a few properties on  the renormalized solutions of
(\ref{1.1})--(\ref{1.3}). Section 3 is devoted to the existence 
result. At last, in Section 4 we prove two uniqueness results. 
The results of the present paper were announced in \cite{BCAG99} and
here we improve the uniqueness results.

\section{Assumptions and definitions}

\subsection{Notations and hypotheses}

In the whole paper, for $q\in \left[ 1,\infty \right[ $ we denote
by  $W_{\Gamma _{d}}^{1,q}( \Omega ) $ the space of
functions belonging to 
$W^{1,q}( \Omega ) $ which have a null trace on $\Gamma
_{d}$. Since $\Omega$ is a bounded and connected open subset of $\R^N$
with Lipschitz boundary and since $\sigma(\Gamma_d)>0$, the space
$W_{\Gamma _{d}}^{1,q}( \Omega ) $ is provided by the norm
$\left\| v\right\| _{W_{\Gamma _{d}}^{1,q}( \Omega 
) }=\left\| Dv\right\| _{L^{q}( \Omega ) }$ (see
e.g. \cite{ziemer}).

\par

We assume that $\aop :\Omega \times \mathbb{R}\longmapsto\mathbb{R}^{N}$ , $ 
\lambda : \Omega \times \mathbb{R}\longmapsto\mathbb{R}$ and $ 
\Phi :\Omega \times \mathbb{R}\longmapsto \mathbb{R}^{N}$ are
Carath\'{e}odory functions such that for $1<p\leq N$ we have :  
\begin{gather}
\exists \alpha >0,\quad \aop ( x,\xi ) \cdot \xi \geq
\alpha \left| \xi \right| ^{p},\text{ }\forall \xi \in \mathbb{R}^{N},\text{\
\ a.e. }x\text{ in }\Omega ;  \label{2.1}
\\
\big( \aop ( x,\xi ) -\aop ( x,\xi ^{\prime
})\big) \cdot ( \xi -\xi ^{\prime })  \geq 0\text{ \ }\forall
\xi,\xi' \in \mathbb{R}^{N},\text{\ \ a.e. }x\text{ in }\Omega ;  \label{2.2}
\\
\left| \aop ( x,\xi ) \right| \leq \beta ( \left|
d( x) \right| +\left| \xi \right| ^{p-1}) \text{ }\forall
\xi \in \mathbb{R}^{N},\text{\ \ a.e. }x\text{ in }\Omega,\text{ with }
d\in L^{p^{\prime }}( \Omega ) ;  \label{2.3}
\\
\lambda ( x,s) s \geq 0\text{ \ }\forall s\in \mathbb{R},\text{\
\ a.e. }x\text{ in }\Omega ;  \label{2.4} \\
\forall k >0,\text{ }\exists c_{k}>0\text{\   such that }\left| \lambda (
x,s) \right| \leq c_{k}\text{ \ }\forall \left| s\right| \leq k,\text{
\ \ a.e. }x\text{ in }\Omega,  \label{2.4bis}
\\
\left| \Phi ( x,s) \right| \leq \left|  b (
x) \right| ( 1+\left| s\right| ) ^{p-1}\text{ \ }\forall
s\in \mathbb{R},\text{\ \ a.e. }x\text{ in }\Omega,  \label{2.5}
\\
\intertext{with }
\left\{ 
\begin{aligned}[c]
 b \in L^{\frac{N}{p-1}}( \Omega ) \text{ if }p<N, \\ 
\int_{\Omega }( 1+\left|  b \right| ) ^{\frac{N}{N-1} 
}( \ln ( 1+\left|  b \right| ) )
^{N-1}dx<\infty \text{ \ if }p=N; 
\end{aligned} 
\right.   \label{2.5bis}
\\
 f \in L^{1}( \Omega ) .  \label{2.5BIS}
\end{gather}
\par
For
any $k\geq 0$, the truncation function at height $\pm k$ is defined by 
\linebreak $ T_{k}( s) =\max ( -k,\min ( s,k) )$.
 For any integer $n\geq 1$, let us define the bounded positive function
\begin{equation}
h_{n}( s) =1-\frac{\left| T_{2n}( s) -T_{n}(
s) \right| }{n}.  \label{2.6}
\end{equation} 
 For
any measurable subset $E$ of $\Omega $, $\nbOne _{E}$ denotes the
characteristic function of the subset $E$.

\subsection{Definition of a renormalized solution of (\ref 
{1.1})--(\ref{1.3})}

Following  \cite{BBGGV} let us recall the definition of the gradient of
functions whose truncates belong to $W_{\Gamma _{d}}^{1,p}( \Omega
)$.

\begin{definition}
\label{def2.1}Let $u$ be a measurable function defined on $\Omega $ which
is finite almost everywhere such that $T_{k}( u) \in W_{\Gamma
_{d}}^{1,p}( \Omega ) $ for every $k>0$. Then there exists a
unique measurable function $v:\Omega \longrightarrow \mathbb{R}^{N}$ such
that 
\begin{equation*}
DT_{k}( u) =\nbOne_{\left\{ \left| u\right| <k\right\} }v\text{ \
\ a.e. \ in }\Omega ,\text{ }\forall k>0.
\end{equation*}
This function $v$ is called the gradient of $u$ and is denoted by $Du$.
\end{definition}

We now give the definition of a renormalized solution
of problem (\ref{1.1})--(\ref{1.3}).

\begin{definition} \label{def2.2} 
A measurable function $u$ defined on $\Omega $ and finite
almost everywhere on $\Omega$  is called a renormalized solution of (\ref 
{1.1})--(\ref{1.3}) if 
\begin{equation}
T_{k}( u) \in W_{\Gamma _{d}}^{1,p}( \Omega ) ,\text{
\ }\forall k>0;  \label{2.7}
\end{equation} 
\begin{equation}
\lim_{n\rightarrow \infty }\frac{1}{n}\int_{\left\{ \left| u\right|
<n\right\} }\aop ( x,Du) \cdot Dudx=0;  \label{2.8}
\end{equation} 
and if for every $h\in W^{1,\infty }( \mathbb{R}) $, with
compact support and any  $\varphi \in W_{\Gamma _{d}}^{1,p}( \Omega
) \cap L^{\infty }( \Omega ) $,  
\begin{multline} \label{2.9}
\int_{\Omega }\lambda ( x,u) \varphi h( u)
dx+\int_{\Omega }h( u) \aop ( x,Du) \cdot
D\varphi dx+\int_{\Omega }h( u) \Phi ( x,u)
\cdot D\varphi dx   \\
+\int_{\Omega }\varphi h^{\prime }( u) \aop (
x,Du) \cdot Dudx+\int_{\Omega }\varphi h^{\prime }( u) 
\Phi ( x,u) \cdot Dudx  
=\int_{\Omega } f \varphi h( u) dx. 
\end{multline}
\end{definition}

\begin{remark}
Condition (\ref{2.7}) and Definition \ref{def2.1} allow to define $Du$
almost everywhere in $\Omega $. Condition (\ref{2.8}) which is crucial to obtain uniqueness
results is standard in the
context of renormalized solution and gives additional information on $Du$
for large value of $|u|$. Equality (\ref{2.9}) is
formally obtained by using in (\ref{1.1}) the test function $\varphi
h( u) $ and taking into account the boundary conditions (\ref{1.2})
and (\ref{1.3}).

Every term in (\ref{2.9}) is well defined. Indeed let $k > 0$ such
that ${\text {supp}}(h)\subset \left[ -k,k\right] $. From assumption
(\ref{2.4bis}) we 
have $\left| \lambda ( x,u) \varphi h( u) \right|
\leq c_{k}\left\| \varphi h\right\| _{L^{\infty }( \Omega ) }$
a.e. in $\Omega $ and then  $\lambda ( x,u) \varphi h(
u) $ lies in $L^{1}( \Omega ) $. Since $h( u) 
(\aop ( x,Du)+ \Phi (
x,u)) =h( u)( \aop (
x,DT_{k}( u) )+ \Phi ( x,T_{k}( u)
))  $ a.e. in $\Omega $, from  (\ref{2.3}), (\ref{2.5}) and (\ref 
{2.7}) it follows that $h( u) (\aop ( x,Du)+
 \Phi ( x,u)) $ belongs to $\big(
L^{p^{\prime }}( \Omega) \big)^{N}$. Thus $h( u) (\aop ( x,Du)+
 \Phi ( x,u))\cdot D\varphi$ is integrable on $\Omega$.
The same arguments imply that  $\varphi h^{\prime }( u) \aop  
( x,Du) \cdot Du=\varphi h^{\prime }( u) \aop  
( x,DT_{k}( u) ) \cdot DT_{k}( u)$
 and $\varphi h^{\prime }( u) 
\Phi ( x,u) \cdot Du=\varphi h^{\prime }( u)  
\Phi ( x,T_{k}( u) ) \cdot DT_{k}(
  u)$ lie  in
$L^{1}( \Omega)$.
At last it is clear that $f \varphi h( u) \in
L^{1}( \Omega )$.
\end{remark}


\subsection{Properties of a renormalized solution of
  (\ref{1.1})--(\ref{1.3}).} 

\bigskip It is  well known (see \cite{DMOP99}) that if $ f \in
L^{1}( \Omega ) $, any renormalized solution of the equation $-\diw  
( \aop ( x,Du) ) = f $  in $\Omega$, 
$u=0$ on $\partial \Omega $ is also a solution in the sense of distribution.
We establish a similar result in Proposition \ref{prop2.6} namely that any function $\psi \in
W_{\Gamma _{d}}^{1,q}(\Omega)$, with $q>N$, is an admissible
test function in (\ref{1.1})--(\ref{1.3}). We first give two technical
lemmas which will be used in the limit case $p=N$.

\begin{lemma}
\label{lem2.4bis}$\forall \omega > 0 $,  $\exists \eta >0$ such 
\begin{equation}
\forall v\in W_{\Gamma _{d}}^{1,N}( \Omega ) ,\quad
\int_{\Omega }\exp \bigg( \Big( \frac{v}{\eta \left\| Dv\right\|
_{(L^{N}( \Omega ))^N}}\Big)^{\frac{N}{N-1}}-1\bigg) dx\leq
\omega.  \label{2.9bis}
\end{equation}
\end{lemma}

\begin{lemma}
\label{lem2.4ter}$\exists c( N)>0 $ such that 
\begin{equation}
\forall \theta >0,\text{ }\forall x,y\in \mathbb{R},\text{ \ }\left|
xy\right| \leq \frac{c( N) }{\theta }\Big( ( 1+\left|
x\right| ) \big( \ln ( 1+\left| x\right| ) \big)
^{N-1}+e^{\left| \theta y\right| ^{\frac{1}{N-1}}}-1\Big).  \label{2.9ter}
\end{equation}
\end{lemma}

\begin{remark}
Property (\ref{2.9bis}) is a consequence of the limit-case of Sobolev's
embedding theorem (see \cite{GT77}). Inequality (\ref{2.9ter}) can be easily derived
by induction using the standard inequality 
$\left| xy\right| \leq ( 1+\left| x\right| ) \ln
( 1+\left| x\right| ) +e^{|y|}-1$, $\forall x,y\in\R$. We
leave the details of the proofs to the reader.
\end{remark}
\par
\smallskip
In the following lemma we give
some regularity results of a renormalized solution of (\ref{1.1})--(\ref 
{1.3}).
\begin{lemma} \label{lem2.4}
Assume that (\ref{2.1})--(\ref{2.5BIS})  hold true. If $u$ is a
renormalized solution of (\ref{1.1})--(\ref{1.3}) then  
\begin{equation}
\lambda ( x,u) \in L^{1}( \Omega ) ;  \label{2.10}
\end{equation} 
\begin{equation}
\lim_{n\rightarrow \infty }\frac{1}{n}\int_{\left\{ \left|
u\right| <n\right\} }\big| \Phi ( x,u)\big| \cdot \big|Du\big|
dx=0;  \label{2.11}
\end{equation} 
\begin{equation}
\left| Du\right| ^{p-1}\in L^{q}( \Omega ) ,\text{ \ \ \ }\forall 
\,1<q<\frac{N}{N-1};  \label{2.12}
\end{equation} 
\begin{gather}
\left| u\right| ^{p-1}\in L^{q}( \Omega ) ,\text{ \ \ \ }\forall 
\,1<q<\frac{N}{N-p};  \label{2.13} \\
\frac{\big|Du\big|^p}{\big(1+|u|\big)^{1+m}} \in L^1(\Omega),
\quad\forall m>0. \label{2.durdur}
\end{gather}
\end{lemma}

\begin{proof}[Sketch of the proof of Lemma  \ref{lem2.4}]

Regularities (\ref{2.12}), (\ref{2.13}) and (\ref{2.durdur}) are easy consequences of the
 estimate techniques of L. Boccardo and T. Gallou\"{e}t developed in
\cite{BG89} (see also \cite{BBGGV}
 and \cite{BMMP99}). Indeed (\ref{2.1}), (\ref{2.7}) and (\ref{2.8}) yield
that  $\int_{\Omega }\left| DT_{k}( u) \right| ^{p}dx\leq ck+L$ $ 
\forall k>0$.
\par
Let us prove (\ref{2.11}). Assumption (\ref{2.5}) and H\"{o}lder's
inequality lead to 
\begin{equation}
\frac{1}{n}\int_{\left\{ \left| u\right| <n\right\} }\big| \Phi  
( x,u)\big| \cdot \big| Du\big| dx\leq \bigg( \frac{1}{n}\int_{\left\{
\left| u\right| <n\right\} }\left|  b \right| ^{p^{\prime }}\big(
1+\left| u\right| \big) ^{p}dx\bigg)^{\frac{1}{p^{\prime }}}
\bigg( \frac{ 
1}{n}\int_{\left\{ \left| u\right| <n\right\} }\left| Du\right|
^{p}dx\bigg)^{\frac{1}{p}}.  \label{2.14}
\end{equation}

If $p<N$, using Sobolev's embedding theorem we have 
\begin{align*}
\bigg( \int_{\left\{ \left| u\right| <n\right\} }\left|  b \right|
^{p^{\prime }}\big( 1+\left| u\right| \big) ^{p}dx \bigg)^{\frac{1}{ 
p^{\prime }}} &\leq  \left\|  b \right\| _{L^{\frac{N}{p-1}}(
\Omega ) }\bigg( \int_{\Omega }\big( 1+\left| T_{n}( u)
\right| \big)^{\frac{pN}{N-p}}dx\bigg)^{\frac{( p-1) (
N-p) }{pN}} \\
 & \leq c\left\|  b \right\| _{L^{\frac{N}{p-1}}( \Omega )
}\Big( 1+\left\| DT_{n}( u) \right\| _{( L^{p}( \Omega
))^{N}}^{p-1}\Big) 
\end{align*} 
and  with (\ref{2.14}) we obtain
\begin{equation}
\frac{1}{n}\int_{\left\{ \left| u\right| <n\right\} }\left| \Phi  
( x,u) \cdot Du\right| dx\leq \frac{c}{n} \left\|  b
\right\| _{L^{\frac{N}{p-1}}(\Omega)} (
1+
\left\| DT_{n}( u) \right\| _{( L^{p}( \Omega
) ) ^{N}}^{p}),   \label{2.15}
\end{equation} 
where $c$ is a constant independent on $n$.
\par

If $p=N$, using Lemma \ref{lem2.4bis} and  Lemma \ref{lem2.4ter}  we obtain
after a few computations, for $\eta >0$, 
\begin{align*}
\null &\int_{\left\{ \left| u\right| <n\right\} }\left|  b \right|
^{N^{\prime }}( 1+\left| u\right| ) ^{N}dx \\
&\leq c\left\|  b \right\| _{L^{\frac{N}{N-1}}( \Omega )
}^{N'}+c\int_{\Omega }\left|  b \right| ^{N^{\prime }}\left|
T_{n}( u) \right| ^{N}dx \\
&\leq c\left\|  b \right\| _{L^{\frac{N}{N-1}}( \Omega )
}^{N'} 
+ cc{( N) }\left\| DT_{n}( u) \right\| _{(
L^{N}(\Omega))^{N}}^{N}\int_{\Omega }\big( 1+\eta
^{N}\left|  b \right| ^{N^{\prime }}\big) \big( \ln
( 1+\eta ^{N}\left|  b \right| ^{N^{\prime }})
\big)^{N-1}dx \\
&\quad +cc{( N) }\left\| DT_{n}( u) \right\| _{(
L^{N}( \Omega ) ) ^{N}}^{N}\int_{\Omega }\bigg( \exp \Big( 
\frac{\left| T_{n}( u) \right| }{\eta \left\| DT_{n}(
u) \right\| _{( L^{N}( \Omega ) ) ^{N}}}
\Big)^{N^{\prime }}-1\bigg) dx,
\end{align*} 
where $c$ and $c(N)$ are positive constants independent of $n$.
Choosing $\omega=1$ in Lemma \ref{lem2.4bis} (then $\eta $ is fixed), we get 
\begin{align}
\int_{\left\{ \left| u\right| <n\right\} }\left|  b \right|
^{N^{\prime }}( 1+\left| u\right| ) ^{N}dx 
\leq {} &
c\left\|  b \right\| _{L^{\frac{N}{N-1}}( \Omega )
}^{N'}+cc{( N) }\left\| DT_{n}( u) \right\| _{(
L^{N}( \Omega ) ) ^{N}}^{N} \notag \\ 
& \quad {} + cc{( N) }\left\| DT_{n}( u) \right\| _{(
L^{N}( \Omega ) ) ^{N}}^{N}\int_{\Omega }( 1+\left| 
 b \right| ) ^{N^{\prime }}( \ln ( 1+\left|  b  
\right| ) ) ^{N-1}dx \notag \\
{} \leq  {} & {} c \big(1 + \|DT_n(u)\|_{(L^N(\Omega))^N}^N\big), \label{durdur1}
\end{align} 
where $c$ is a positive constant depending on $N$, $\eta$, $b$ and $\Omega $.
\par
From (\ref{2.14}) together with (\ref{2.15}) and (\ref{durdur1}) it
follows that in both cases $p<N$ and $N=p$
\begin{equation}\label{durdur2}
\frac{1}{n}\int_{\left\{ \left| u\right| <n\right\} }\left| \Phi  
( x,u) \cdot Du\right| dx\leq \frac{c}{n}( 
1
+\left\| DT_{n}( u) \right\| _{( L^{p}( \Omega
) ) ^{N}}^{p}) \quad\forall n\geq 1,
\end{equation} 
where $c$ is a positive constant which depends on $N$, $p$, $b$ and
$\Omega$.
 Assumption  (\ref{2.1}), condition (\ref{2.8}) and (\ref{durdur2}) 
lead to (\ref{2.11}).
\par

We are now in a position to obtain (\ref{2.10}). For any $n>0$, the function $h_{n}$
(see (\ref{2.6})) is Lipschitz continuous with compact support, so
that (\ref{2.9}) yields with $\varphi =T_{1}( u) \in W_{\Gamma 
_{d}}^{1,p}( \Omega ) \cap L^{\infty }( \Omega)$ 
\begin{multline*}
\int_{\Omega }\lambda ( x,u) T_{1}( u) h_{n}(
u) dx+\int_{\Omega }h_{n}( u) \aop ( x,Du)
\cdot DT_{1}( u) dx \\
+\int_{\Omega }h_{n}( u) \Phi ( x,u) \cdot
DT_{1}( u) dx+\int_{\Omega }T_{1}( u) h_{n}^{\prime
}( u) \aop ( x,Du) \cdot Dudx \\
+\int_{\Omega }T_{1}( u) h_{n}^{\prime }( u)  
\Phi ( x,u) \cdot Dudx 
=\int_{\Omega } f T_{1}( u) h_{n}( u) dx.
\end{multline*} 
Because $0\leq h_{n}( u) \leq 1$ and 
$\aop ( x,Du) \cdot DT_{1}( u) \geq 0$ almost
everywhere in $\Omega $, we deduce that
\begin{multline*}
\int_{\Omega }\lambda ( x,u) T_{1}( u) h_{n}(
u) dx  \label{2.18} 
\leq \left\|  f \right\| _{L^{1}( \Omega ) }
+\int_{\Omega }\left| \Phi ( x,u)
\cdot DT_{1}( u) \right| dx 
\\
+\frac{1}{n} 
\int_{\left\{ \left| u\right| <n\right\} }( \aop (
x,Du) \cdot Du+\left| \Phi ( x,u) \cdot 
Du\right| ) dx.
\end{multline*} 
Since $u$ is finite almost everywhere in $\Omega $,  $h_{n}( u) $
converges to 1 almost everywhere in $\Omega $ and it is bounded by 1.
Assumption (\ref{2.4}) and Fatou lemma together with (\ref{2.8}) and (\ref 
{2.11}) allow to conclude that 
\begin{equation}
\int_{\Omega }\lambda ( x,u) T_{1}( u) dx\leq \left\| 
 f \right\| _{L^{1}( \Omega ) }+\int_{\Omega }\left| 
\Phi ( x,u) \cdot DT_{1}( u) \right| dx.
\label{2.19}
\end{equation} 
At last writing $\left| \lambda ( x,u) \right| \leq $ $\lambda
( x,u) T_{1}( u) +\left| \lambda ( x,u)
\right| \indi{\left\{ \left| u\right| \leq 1\right\} }$ a.e. in $ 
\Omega $, and using (\ref{2.4bis}), (\ref{2.15}) and (\ref{2.19}) lead
to (\ref{2.10}).
\end{proof}

\begin{proposition} \label{prop2.6} 
Assume that (\ref{2.1})--(\ref{2.5BIS}) hold true. If $u
$ is a renormalized solution of (\ref{1.1})--(\ref{1.3}) then for any $\psi
\in \underset{r>N}{\bigcup }W_{\Gamma _{d}}^{1,r}( \Omega ) $ we
have 
\begin{equation}\label{durdurprop2.6}
\int_{\Omega }\lambda ( x,u) \psi dx+\int_{\Omega }\aop  
( x,Du) \cdot D\psi dx+\int_{\Omega }\Phi (
x,u) \cdot D\psi dx=\int_{\Omega } f \psi dx.
\end{equation}
\end{proposition}

\begin{proof}[Sketch of proof]
Let $\psi \in W_{\Gamma _{d}}^{1,r}( \Omega ) $ with $r>N$.
Sobolev's embedding theorem implies that $\varphi \in W_{\Gamma
_{d}}^{1,r}( \Omega ) \cap L^{\infty }( \Omega ) $ and
(\ref{2.9}) with $h=h_{n}$ leads to 
\begin{multline}\label{durdur3}
\int_{\Omega }\lambda ( x,u) \psi h_{n}( u)
dx+\int_{\Omega }h_{n}( u) \aop ( x,Du) \cdot
D\psi dx 
+\int_{\Omega }h_{n}( u) \Phi ( x,u) \cdot
D\psi dx
\\
+\int_{\Omega }\psi h_{n}^{\prime }( u) \aop (
x,Du) \cdot Dudx 
+\int_{\Omega }\psi h_{n}^{\prime }( u) \Phi (
x,u) \cdot Dudx=\int_{\Omega } f \psi h_{n}( u) dx.
\end{multline} 

Assumptions (\ref{2.3}) and (\ref{2.12}) lead to 
\begin{equation*}
\aop ( x,Du) \in ( L^{q}( \Omega ) )
^{N}\text{ \ \ }\forall 1\leq q<\frac{N}{N-1},
\end{equation*} 
and then $\aop ( x,Du) \cdot D\psi \in L^{1}( \Omega
)$. Similarly (\ref{2.5}) and (\ref{2.13}) give that $\Phi  
( x,u) \cdot D\psi \in L^{1}( \Omega )$. 
Since $h_n(u)$ converges to 1 almost everywhere in $\Omega$ and it is
bounded by 1 and recalling that $|h'_n(s)|=1/n\indi{\{n<|s|<2n\}}(s)$ a.e. on $\R$, it
is then a straightforward task to pass the limit in (\ref{durdur3})
using Lebesgue's dominated convergence theorem, (\ref{2.8}),
(\ref{2.10}) and (\ref{2.11}). Such a limit process leads to  (\ref{durdurprop2.6}).
\end{proof}

\section{Existence of a renormalized solution}

\begin{theorem}
\bigskip \label{th3.1}Under assumptions (\ref{2.1})--(\ref{2.5BIS})  there
exists a renormalized solution of equation (\ref{1.1})--(\ref{1.3}).
\end{theorem}

\begin{proof}
The proof relies on passing to the limit in an approximate problem.
\par
{\slshape Step 1.} For $\varepsilon >0$, let us define
\begin{equation}
\ \lambda _{\varepsilon }( x,s) =\varepsilon \left| s\right|
^{p-2}s+\lambda \big( x,T_{{1/\varepsilon }}( s) \big),
\label{3.1}
\end{equation} 
\begin{equation}
\Phi _{\varepsilon }( x,s) =\Phi \big( x,T_{ 
{1/\varepsilon }}(s)\big)   \label{3.2}
\end{equation} 
and $f^{\varepsilon }\in L^{p^{\prime }}( \Omega ) $
such that 
\begin{equation}
 f ^{\varepsilon }\overset{\varepsilon \rightarrow 0}{\longrightarrow 
} f \qquad \text{strongly in }L^{1}( \Omega ).   \label{3.3}
\end{equation} 
From the classical results of Leray--Lions \cite{Lio69}, an application of
the Leray-Schauder  fixed point theorem allows to show that for any $ 
\varepsilon >0$ there exists $u^{\varepsilon }\in W_{\Gamma
_{d}}^{1,p}( \Omega ) $ such that $\forall \psi \in
W_{\Gamma _{d}}^{1,p}( \Omega )$ 
\begin{equation}
\int_{\Omega }\lambda _{\varepsilon }( x,u^{\varepsilon }) \psi
dx+\int_{\Omega }\aop ( x,Du^{\varepsilon }) \cdot D\psi
dx+\int_{\Omega }\Phi _{\varepsilon }( x,u^{\varepsilon
}) \cdot D\psi dx=\int_{\Omega } f ^{\varepsilon }\psi dx.
\label{3.4}
\end{equation}

We now derive some estimates on $u^{\varepsilon }$. Using $\psi =T_{k}(
u^{\varepsilon }) $ for $k>0$ in (\ref{3.4}) we obtain 
\begin{multline*}
\int_{\Omega }\lambda _{\varepsilon }( x,u^{\varepsilon })
T_{k}( u^{\varepsilon }) dx+\int_{\Omega }\aop (
x,DT_{k}( u^{\varepsilon }) ) \cdot DT_{k}(
u^{\varepsilon }) dx 
\\
+ \int_{\Omega }\Phi _{\varepsilon }(
x,T_{k}( u^{\varepsilon }) ) \cdot DT_{k}(
u^{\varepsilon }) dx=\int_{\Omega } f ^{\varepsilon
}T_{k}( u^{\varepsilon }) dx.
\end{multline*} 
From the coercivity of the operator $\aop$, the positivity of the
function $\lambda $ and (\ref{2.5}) it follows that 
\begin{multline*}
\int_{\Omega }\lambda ( x,u^{\varepsilon }) T_{k}(
u^{\varepsilon }) dx+\alpha \int_{\Omega }\left| DT_{k}(
u^{\varepsilon }) \right| ^{p}dx
\\
\leq \int_{\Omega } f  
^{\varepsilon }T_{k}( u^{\varepsilon }) dx+\int_{\Omega } b 
( x) ( 1+\left| T_{k}( u^{\varepsilon }) \right|
) ^{p-1}\left| DT_{k}( u^{\varepsilon }) \right| dx
\end{multline*} 
and then Young's inequality gives 
\begin{equation}
\int_{\Omega }\lambda ( x,u^{\varepsilon }) T_{k}(
u^{\varepsilon }) dx+\int_{\Omega }\left| DT_{k}( u^{\varepsilon
}) \right| ^{p}dx\leq ( M+1) ( k+k^{p}) 
\label{3.5}
\end{equation} 
where $M$ is a generic constant independent of $k$ and $\varepsilon$.
Inequality (\ref{3.5}) implies that $\forall k>0$
\begin{equation}
T_{k}( u^{\varepsilon }) \text{ is bounded in } 
W_{\Gamma _{d}}^{1,p}( \Omega ) \label{3.6}
\end{equation} 
and due to (\ref{2.4}) and (\ref{2.4bis}) 
\begin{equation}
\lambda ( x,u^{\varepsilon }) \text{ \ \ is bounded in } 
L^{1}( \Omega ) .\label{3.7}
\end{equation}

As a consequence of (\ref{2.3}), (\ref{3.6}) and (\ref{3.7}) there exists a
subsequence  (still denoted by $\varepsilon $) and a measurable function  
$u\,:\,\Omega \longrightarrow \overline{\mathbb{R}}$ such that
$\lambda(x,u)\in L^1(\Omega)$ and such that
\begin{gather}
u^{\varepsilon }\overset{\varepsilon \rightarrow 0}{\longrightarrow }u\text{
\ almost everywhere in }\Omega,\label{3.8} 
\\
\forall k >0\qquad T_{k}( u^{\varepsilon }) \overset{ 
\varepsilon \rightarrow 0}{\rightharpoonup }T_{k}( u) \text{
weakly in }W_{\Gamma _{d}}^{1,p}( \Omega ) ,  \label{3.9} 
\\
\forall k >0,\quad
\aop ( x,DT_{k}(
u^{\varepsilon }) ) \overset{\varepsilon \rightarrow 0}{ 
\rightharpoonup }\sigma_{k}\text{ weakly in }( L^{p^{\prime
}}( \Omega ) ) ^{N},  \label{3.9bis}
\end{gather}
where $ \sigma_{k}\in \big( L^{p^{\prime }}(
\Omega ) \big) ^{N}$.

We claim that $u$ is finite almost everywhere in $\Omega $ (remark that if 
$\lambda ( x,s) =\lambda \left| s\right| ^{p-2}s$, with $\lambda
>0$,  it is obvious since $\lambda ( x,u) \in L^{1}( \Omega
)$) through a ``log--type'' estimate on $u^\varepsilon$ (such a ``log--type'' estimate
is also performed in \cite{B}, see also \cite {BOP,droniou02,droniou03}). 
Let us consider the real valued function 
\begin{equation*}
\psi _{p}( r)
=\int_{0}^{r}\frac{ds}{( 1+\left| s\right| ) ^{p}}\quad\forall r\in \mathbb{R}.
\end{equation*} 
Since $p>1$, $\psi _{p}( u) \in W_{\Gamma _{d}}^{1,p}(
\Omega ) \cap L^{\infty }( \Omega ) $ with $\left\| \psi
_{p}( u) \right\| _{L^{\infty }( \Omega ) }\leq \frac{1 
}{p-1}$ and  $D\psi _{p}( u) =\frac{Du}{( 1+\left|
u\right| ) ^{p}}$ almost everywhere in $\Omega$, $\psi
_{p}( u^{\varepsilon }) $ is an admissible test function in (\ref 
{3.4}). It follows that
\begin{equation*}
\int_{\Omega }\lambda _{\varepsilon }( x,u^{\varepsilon }) \psi
_{p}( u^{\varepsilon }) dx+\int_{\Omega }\aop (
x,u^{\varepsilon }) \cdot D\psi _{p}( u^{\varepsilon })
dx+\int_{\Omega }\Phi _{\varepsilon }( x,u^{\varepsilon
}) \cdot D\psi _{p}( u^{\varepsilon }) dx\leq \frac{M}{p-1}
\end{equation*} 
and due to the definition of $\lambda _{\varepsilon }$ together with  (\ref{2.1}),
(\ref{2.4}) and (\ref{2.5}) we have 
\begin{equation*}
\alpha \int_{\Omega }\frac{\left| Du^{\varepsilon }\right| ^{p}}{(
1+\left| u^{\varepsilon }\right| ) ^{p}}dx\leq \frac{M}{p-1} 
+\int_{\Omega } b ( x) \frac{\left| Du^{\varepsilon
}\right| }{( 1+\left| u^{\varepsilon }\right| ) }dx
\end{equation*} 
where $M$ is a generic constant independent of $\varepsilon$. Young's
inequality leads to
\begin{equation*}
\int_{\Omega }\frac{\left| Du^{\varepsilon }\right| ^{p}}{( 1+\left|
u^{\varepsilon }\right| ) ^{p}}dx\leq M\Big( 1+\int_{\Omega }\left| 
 b ( x) \right| ^{\frac{p}{p-1}}dx\Big) .
\end{equation*} 
Since $p\leq N$ and $u^{\varepsilon }\in W_{\Gamma _{d}}^{1,p}( \Omega
)$, the regularity of the function $ b $ (see (\ref{2.5bis}))
implies that the field $\ln ( 1+\left| u^{\varepsilon }\right| ) $
is bounded in $W_{\Gamma _{d}}^{1,p}( \Omega ) $ uniformly with
respect to $\varepsilon$. From (\ref{3.8}) and (\ref{3.9}) it follows that $ 
\ln ( 1+\left| u\right| ) $ belongs to $W_{\Gamma
_{d}}^{1,p}( \Omega ) $ and then $u$ is finite almost everywhere
in $\Omega$.
\par
\medskip
{\slshape Step 2.} We prove the following lemma.

\begin{lemma}
\label{lemme3.2} 
\begin{equation}
{\lim_{n\rightarrow \infty }}{\limsup_{\varepsilon
\rightarrow 0}}\frac{1}{n}\int_{\left\{ \left| u^{\varepsilon }\right|
<n\right\} }\left| Du^{\varepsilon }\right| ^{p}dx=0.  \label{3.10}
\end{equation}
\end{lemma}

\begin{proof}
Taking the admissible test function $T_{n}(u^{\varepsilon
})/n$ (which belongs to $W_{\Gamma _{d}}^{1,p}( \Omega )
\cap L^{\infty }( \Omega ) $) in (\ref{3.4}) yields that
\begin{multline*}
\frac{1}{n}\int_{\Omega }\lambda _{\varepsilon }(
x,u^{\varepsilon }) T_{n}( u^{\varepsilon })
dx+
\frac{1}{n}\int_{\Omega }\aop ( x, Du^{\varepsilon }
) \cdot DT_{n}( u^{\varepsilon }) dx 
\\
+\frac{1}{n}
\int_{\Omega } 
\Phi_{\varepsilon }( x, u^{\varepsilon } )
\cdot DT_{n}( u^{\varepsilon }) dx 
=\frac{1}{n} 
\int_{\Omega } f ^{\varepsilon }T_{n}( u^{\varepsilon })
dx.
\end{multline*}
Since $\lambda _{\varepsilon }( x,u^{\varepsilon }) T_{n}(
u^{\varepsilon }) \geq 0$ almost everywhere in $\Omega$, using
(\ref{2.1})  and (\ref{2.5}) we get 
\begin{multline} \label{3.11}
\frac{\alpha }{n}\int_{\Omega }\left| DT_{n}( u^{\varepsilon })
\right| ^{p}dx\leq \frac{1}{n}\Big( \int_{\Omega }\left|  b (
x) \right| ( 1+\left| T_{n}( u^{\varepsilon }) \right|
) ^{p-1}\left| DT_{n}( u^{\varepsilon }) \right|
dx \\
{} + \int_{\Omega } f ^{\varepsilon }T_{n}( u^{\varepsilon
}) dx\Big).
\end{multline} 
As a consequence of (\ref{3.3}) and (\ref{3.8}) and using the fact
that $u$ is finite almost
everywhere in $\Omega$,  Lebesgue's convergence theorem leads to
\begin{equation*}
\lim_{n\rightarrow \infty }\limsup_{\varepsilon
\rightarrow 0}\frac{1}{n}\int_{\Omega }\left|  f ^{\varepsilon
}T_{n}( u^{\varepsilon }) \right| dx=\lim_{n\rightarrow
\infty }\frac{1}{n}\int_{\Omega }\left|  f T_{n}( u)
\right| dx=0.
\end{equation*} 
Due to (\ref{3.2}), (\ref{3.8}) and (\ref{3.9})   we have 
\begin{equation}
\limsup_{\varepsilon \rightarrow 0}\frac{1}{n}\int_{\Omega }\left| 
 b ( x) \right| ( 1+\left| T_{n}( u^{\varepsilon
}) \right| ) ^{p-1}\left| DT_{n}( u^{\varepsilon })
\right| dx=\frac{1}{n}\int_{\Omega }\left|  b ( x) \right|
( 1+\left| T_{n}( u) \right| ) ^{p-1}\left|
DT_{n}( u) \right| dx.  \label{3.12}
\end{equation}  
We will prove in the sequel by splitting techniques that 
\begin{equation}
\forall \eta >0,\qquad \frac{1}{n}\int_{\Omega }\left|  b (
x) \right| ( 1+\left| T_{n}( u) \right| )
^{p-1}\left| DT_{n}( u) \right| dx\leq \omega _{\eta }(
n) +\frac{\eta }{n}\int_{\Omega }\left| DT_{n}( u) \right|
^{p}dx,  \label{3.13}
\end{equation} 
where $\omega _{\eta }( n) $ goes to zero as $n$ goes to
infinity. Since $T_{n}( u^{\varepsilon }) $ converges to $ 
T_{n}( u) $ weakly in $W_{\Gamma _{d}}^{1,p}( \Omega ) 
$, choosing $\eta $ small enough in (\ref{3.13}) together with (\ref{3.11})
give (\ref{3.10}).

In order to complete the proof, it remains to prove 
(\ref{3.13}). 
Let $\eta >0$ and $R>0$ ($R$ will be fixed later).  By
denoting $E_{R}$ the measurable set $E_{R}=\{ x\in \Omega \,;\,
|u(x)| >R\} $, we write
\begin{equation} \label{3.14}
\frac{1}{n}\int_{\Omega }\left|  b ( x) \right| (
1+\left| T_{n}( u) \right| ) ^{p-1}\left| DT_{n}(
u) \right| dx = I_1(n,R)+ I_2(n,R)
\end{equation}
with
\begin{gather*}
I_1(n,R)= \frac{1}{n}\int_{\Omega
\setminus E_{R}}\left|  b ( x) \right| ( 1+\left|
  T_{n}( u) \right| ) ^{p-1} 
 \left| DT_{R}(
u) \right| dx 
\\
\intertext{and}
I_2(n,R)= \frac{1}{n}\int_{E_{R}}\left|  b ( x)
\right| ( 1+\left| T_{n}( u) \right| ) ^{p-1}\left|
DT_{n}( u) \right| dx. 
\end{gather*} 
Since $T_R(u)\in W^{1,p}_{\Gamma_d}(\Omega)$ we have
$\lim_{n\rightarrow \infty } I_1(n,R)=0$, $\forall R>0$.
\par
We now deals with $I_2(n,R)$ by
distinguishing  the cases $p<N$ and $p=N$.

{\noindent \slshape First case.} Assuming that $p<N$, H\"{o}lder's inequality and Sobolev's embedding
theorem lead to, $\forall n\geq 1$,
\begin{align}
I_2(n,R) &\leq \frac{1}{n}\left\|  
b\right\| _{L^{\frac{N}{p-1}}( E_{R}) }\big\| 1+|T_{n}(
u)| \big\| _{L^{\frac{Np}{N-p}}( \Omega ) }^{p-1}\left\|
DT_{n}( u) \right\| _{( L^{p}( \Omega ) )
^{N}}  \notag \\
&\leq \frac{M}{n}\left\|  b \right\| _{L^{\frac{N}{p-1}}(
E_{R}) }\big( 1+\left\| DT_{n}( u)
\right\| _{( L^{p}( \Omega ) ) ^{N}}^{p}\big) ,
\label{3.15}
\end{align} 
where $M$ is a generic constant not depending on $n$ and $R$. Since $u$ is
finite almost everywhere in $\Omega $ and since $b\in L^{\frac{N}{p-1} 
}( \Omega )$, let  $R>0$ such that $M\|b\|
_{L^{\frac{N}{p-1}}( E_{R}) }<\eta$. Due to 
(\ref{3.14}) and (\ref{3.15}) we obtain (\ref{3.13}).

\par\smallskip

{\noindent \slshape Second case.} We assume that $p=N$. Let us define  $A_{n}=\left\| DT_{n}(
u) \right\| _{( L^{N}( \Omega ) ) ^{N}}$ and let 
$\rho >0$ ($\rho $ will be fixed in the sequel). A few computations and
Lemma \ref{lem2.4ter} (with $\theta =A_{n}^{-N}$) give, $\forall n \geq 1$,
\begin{align}
I_2(n,R) & \leq \frac{2^{N-2}}{n}\left[
\int_{E_{R}}\left|  b ( x) \right| \left| DT_{n}(
u) \right| dx+\int_{E_{R}}\left|  b ( x) \right|
\left| T_{n}( u) \right| ^{p-1}\left| DT_{n}( u)
\right| dx\right]   
\notag 
\\
& \leq \frac{2^{N-2}}{n} A_n 
\Bigg[ \left\|  b \right\| _{L^{\frac{N}{N-1}%
}( \Omega ) }+\bigg( \int_{E_{R}} \big( \rho ^{N-1}\left|  
b( x) \right| \big) ^{\frac{N}{N-1}}\Big( \frac{\left|
T_{n}( u) \right| }{\rho }\Big)^{N} dx \bigg) ^{{(N-1)}/{N}} 
\Bigg]
\notag 
\\
& \leq \frac{2^{N-2}}{n} A_n
\Bigg[ \|b\| _{L^{N^{\prime
}}( \Omega)} 
+ \Bigg( A_{n}^{N}C( N) \bigg[ 
\int_{E_{R}}\bigg\{ \exp \Big[ \Big( \frac{\left| T_{n}( u)
\right| }{\rho A_{n}}\Big) ^{\frac{N}{N-1}}-1\Big] \bigg\} dx 
\notag 
\\
&\qquad\qquad\qquad {} + 
\int_{E_{R}} \Big( 1+\rho ^{N}\left|  b ( x) \right| ^{ 
\frac{N}{N-1}}\Big) \ln \Big( 1+\rho ^{N}\left|  b ( x)
\right| ^{\frac{N}{N-1}}\Big)^{N-1} dx 
\bigg] \Bigg) ^{{(N-1)}/{N}} \Bigg] 
\notag 
\\
 & \leq \frac{2^{N-2}}{n} A_n
\|b\| _{L^{N^{\prime
}}( \Omega)} +\frac{2^{N-2}C( N) }{n} A_n^N \bigg[ 
\int_{\Omega}\bigg\{ \exp \Big[ \Big( \frac{\left| T_{n}( u)
\right| }{\rho A_{n}}\Big) ^{\frac{N}{N-1}}-1\Big] \bigg\} dx \bigg]^{{(N-1)}/{N}}
\notag \\
& \quad {} + \frac{2^{N-2}C( N) }{n}M( \rho
  ,N) A_n^N \bigg[\int_{E_{R}}\big( 1+\left|  b ( x) 
 \right| \big)^{\frac{N}{N-1}}\ln \big( 1+\left|  b (
 x) \right| \big) ^{N-1}dx \bigg]^{{(N-1)}/{N}},
\label{3.16} 
\end{align} 
where $C(N) >0$ is a constant only depending on $N$ (from 
Lemma \ref{lem2.4ter}) and $M(\rho,N)$ only depends on $\rho$ and $N$.
Using Lemma \ref{lem2.4bis} and since $T_{n}( u) $
lies in 
$W_{\Gamma _{d}}^{1,N}(\Omega)$ we can choose firstly $\rho>0$
such that the quantity 
\begin{equation*}
 2^{N-2}C( N)\bigg[ \int_{\Omega}\bigg\{ \exp \Big[ \Big( \frac{\left| T_{n}( u)
\right| }{\rho A_{n}}\Big) ^{\frac{N}{N-1}}-1\Big] \bigg\} dx \bigg]^{\frac{N-1}{N}}
\end{equation*}
 is small enough independently of $n$. Secondly since $u$ is
finite almost everywhere in $\Omega$ (i.e. $\lim_{R\rightarrow \infty
}\text{meas}( E_{R}) =0$) we can choose $R>0$ such that the
quantity
\begin{equation*}
 2^{N-2}C( N) M( \rho ,N) \bigg[\int_{E_{R}}\big( 1+\left|  b ( x)
\right| \big) ^{\frac{N}{N-1}}\ln \big( 1+\left|  b (
x) \right| \big) ^{N-1}dx\bigg]^{\frac{N-1}{N}}
\end{equation*}
is small enough (notice that it is
crucial to choose $\rho $ before choosing $R$). At last we deduce from
(\ref{3.16}) that there exists $R>0$ such that  $\forall n\geq 1$
\begin{align*}
I_2(n,R)
&\leq \frac{2^{N-2}}{n}\left\| 
 b \right\| _{L^{N^{\prime }}( \Omega ) }\left\|
DT_{n}( u) \right\| _{( L^{N}( \Omega ) )
^{N}}+\frac{\eta }{2n}\left\| DT_{n}( u) \right\| _{(
L^{N}( \Omega ) ) ^{N}}^{N} \notag \\
&\leq \frac{ 2^{N+1}\left\|  b \right\|_{L^{N^{\prime
}}( \Omega ) }^{N'}}{n\eta^{1/(N-1)}}+\frac{\eta }{n} 
\left\| DT_{n}( u) \right\| _{( L^{N}( \Omega )
) ^{N}}^{N}.
\end{align*} 
From (\ref{3.14}) and the behavior
of $I_1(n,R)$ as $n$ goes to infinity it follows that (\ref{3.13})
holds true.
\end{proof}

{\slshape Step 3.} We are now in a position to prove the
following lemma.

\begin{lemma} \label{lemme3.3} 
For any $k>0$,  
\begin{equation*}
\lim_{\varepsilon \rightarrow 0}\int_{\Omega } \big( \aop  
( x,DT_{k}( u^{\varepsilon }) ) -\aop (
x,DT_{k}( u) ) \big)\cdot \big( DT_{k}( u^{\varepsilon
}) -DT_{k}( u) \big) dx =0.
\end{equation*}
\end{lemma}

\begin{proof}[Proof of Lemma \ref{lemme3.3}]
The proof relies on similar techniques developed in \cite{BL93}. 
\par
 Let $k$ be a
positive real number and let $n\geq 1$. Using the test function $(
T_{k}( u^{\varepsilon }) -T_{k}( u) )
h_{n}( u^{\varepsilon }) $ which belongs to $W_{\Gamma
_{d}}^{1,p}( \Omega ) \cap L^{\infty }( \Omega ) $
yields that 
\begin{multline}\label{3.18}
\int_{\Omega }\lambda _{\varepsilon }( x,u^{\varepsilon })
( T_{k}( u^{\varepsilon }) -T_{k}( u) )
h_{n}( u^{\varepsilon }) dx 
\\
+\int_{\Omega }h_{n}(
u^{\varepsilon }) \aop ( x,D u^{\varepsilon
}) \cdot D( T_{k}( u^{\varepsilon })
-T_{k}( u) ) dx  
 +\int_{\Omega }( T_{k}( u^{\varepsilon }) -T_{k}(
u) ) h_{n}^{\prime }( u^{\varepsilon }) \aop  
( x,Du^{\varepsilon }) \cdot Du^{\varepsilon
}dx
\\
+\int_{\Omega }h_{n}( u^{\varepsilon }) \Phi_{\varepsilon
}( x,u^{\varepsilon}) \cdot 
D( T_{k}( u^{\varepsilon }) -T_{k}( u) ) dx
+\int_{\Omega }( T_{k}( u^{\varepsilon }) -T_{k}(
u) ) h_{n}^{\prime }( u^{\varepsilon }) \Phi  
_{\varepsilon }( x,u^{\varepsilon } ) \cdot
Du^{\varepsilon }dx 
\\
=\int_{\Omega } f ^{\varepsilon }( T_{k}( u^{\varepsilon
}) -T_{k}( u) ) h_{n}( u^{\varepsilon }) dx
\end{multline} 
We study in the sequel the behavior of each term of (\ref{3.18})
as $\varepsilon \rightarrow 0$ and $n\rightarrow \infty$.\par
Since $h_{n}$ has
a compact support, condition (\ref{2.4bis}) implies that the field $ 
h_{n}( u^{\varepsilon }) \lambda _{\varepsilon }(
x,u^{\varepsilon }) $ is bounded in $L^{\infty }( \Omega ) .
$ Moreover $T_{k}( u^{\varepsilon }) -T_{k}( u) $
converges to $0$ almost everywhere in $\Omega$ and in $L^{\infty }(
\Omega ) $ weak--$*$ as $\varepsilon$ goes to zero. Therefore  we obtain 
\begin{equation}
\lim_{\varepsilon \rightarrow 0}\int_{\Omega }\lambda _{\varepsilon }(
x,u^{\varepsilon }) ( T_{k}( u^{\varepsilon })
-T_{k}( u) ) h_{n}( u^{\varepsilon }) dx=0,
\label{3.19}
\end{equation} 
and similarly one has 
\begin{equation*}
\lim_{\varepsilon \rightarrow 0}\int_{\Omega } f ^{\varepsilon
}( T_{k}( u^{\varepsilon }) -T_{k}( u) )
h_{n}( u^{\varepsilon }) dx=0.
\end{equation*} 
Recalling that $h'_n(r)=-\indi{\{n<|r|<2n\}} \text{sign}(r)/n$ a.e. on $\R$,
from assumption (\ref{2.3}) it follows that 
\begin{multline*}
\Bigg|\int_{\Omega }( T_{k}( u^{\varepsilon }) -T_{k}(
u) ) h_{n}^{\prime }( u^{\varepsilon }) \aop  
( x,Du^{\varepsilon }) \cdot Du^{\varepsilon
}dx\Bigg|
\\
\leq \frac{M}{n}\Big( \int_{\left\{ n<\left| u^{\varepsilon }\right|
<2n\right\} }\left| Du^{\varepsilon }\right| ^{p}dx+\int_{\Omega }\left| 
d( x) \right| ^{p^{\prime }}dx\Big) ,
\end{multline*} 
with $M>0$ not depending on $\varepsilon $ and $n$. Lemma \ref{lemme3.2} and
the regularity of $d$ allow us to conclude that 
\begin{equation}
\lim_{n\rightarrow \infty }\limsup_{\varepsilon
\rightarrow 0}\Bigg| \int_{\Omega }( T_{k}( u^{\varepsilon })
-T_{k}( u) ) h_{n}^{\prime }( u^{\varepsilon }) 
\aop ( x,DT_{n}( u^{\varepsilon }) ) \cdot
Du^{\varepsilon }dx\Bigg| =0.  \label{3.21}
\end{equation} 
In view of (\ref{2.5}) and since $h_{n}$ has a compact support we get $ 
\left| h_{n}( u^{\varepsilon }) \Phi _{\varepsilon
}( x,u^{\varepsilon }) \right| \leq Mb(x)$ almost everywhere in $ 
\Omega$ where $M$ is a constant independent of $\varepsilon$. Moreover the
pointwise convergence of $u^{\varepsilon }$ and the definition of
$\Phi_\varepsilon$ give that 
\begin{equation*}
h_{n}( u^{\varepsilon }) \Phi _{\varepsilon }(
x,u^{\varepsilon }) \overset{\varepsilon \rightarrow 0}{ 
\longrightarrow }h_{n}( u) \Phi ( x,u) \text{
a.e. in }\Omega .
\end{equation*} 
Thus the regularity of $ b $ and Lebesgue's convergence theorem imply that $ 
h_{n}( u^{\varepsilon }) \Phi _{\varepsilon }(
x,u^{\varepsilon }) $ converges to $h_{n}( u)  
\Phi ( x,u) $ strongly in $L^{\frac{p}{p-1}}( \Omega
)$. Due to (\ref{3.9}) we conclude that 
\begin{equation}
\lim_{\varepsilon \rightarrow 0}\int_{\Omega }h_{n}( u^{\varepsilon
}) \Phi _{\varepsilon }( x, u^{\varepsilon
}) \cdot D( T_{k}( u^{\varepsilon })
-T_{k}( u) ) dx=0  \label{3.22}
\end{equation} 
and similar arguments lead to 
\begin{equation}
\lim_{\varepsilon \rightarrow 0}\int_{\Omega }( T_{k}(
u^{\varepsilon }) -T_{k}( u) ) h_{n}^{\prime }(
u^{\varepsilon }) \Phi _{\varepsilon }( x, u^{\varepsilon
}) \cdot Du^{\varepsilon }dx=0.  \label{3.23}
\end{equation} 
From  (\ref{3.18}) together with (\ref{3.19})--(\ref{3.23}) it follows
that 
\begin{equation}
\lim_{n\rightarrow \infty }{\limsup_{\varepsilon
\rightarrow 0}}\int_{\Omega }h_{n}( u^{\varepsilon }) \aop  
( x,D u^{\varepsilon } ) \cdot D(
T_{k}( u^{\varepsilon }) -T_{k}( u) ) dx=0.
\label{3.24}
\end{equation} 
\par
Since $\aop ( x,0) =0$ almost everywhere in $\Omega $ we
have for $k^{\prime }>k$ 
\begin{equation*}
\aop ( x,DT_{k}( u^{\varepsilon }) ) =\nbOne_{\left\{ \left|
      u^{\varepsilon }\right| <k\right\} }\aop ( 
x,DT_{k^{\prime }}( u^{\varepsilon }) ) \quad\text{almost everywhere
  in $\Omega$}
\end{equation*} 
and due to (\ref{3.8}) and (\ref{3.9bis}) we get 
\begin{equation*}
\sigma_{k}=\nbOne_{\left\{ \left| u\right| <k\right\} } 
\sigma_{k^{\prime }}\text{ a.e. on }\Omega \setminus \left\{ \left|
u\right| =k\right\} .
\end{equation*} 
Since $DT_{k}( u) =0$ a.e. on $\left\{ \left| u\right| =k\right\} 
$ we obtain that 
\begin{equation}
\sigma_{k}\cdot DT_{k}( u) =\sigma_{k^{\prime
}}\cdot DT_{k}( u) \text{ a.e. on }\Omega ,\label{3.25}
\end{equation} 
and then if $n\geq k$ 
\begin{align*}
\lim_{\varepsilon \rightarrow 0}\int_{\Omega }h_{n}( u^{\varepsilon
}) \aop ( x,Du^{\varepsilon }) \cdot DT_{k}(
u) dx =&\lim_{\varepsilon \rightarrow 0}\int_{\Omega }h_{n}(
u^{\varepsilon }) \aop ( x,DT_{2n}( u^{\varepsilon
}) ) \cdot DT_{k}( u) dx  \notag \\
=&\int_{\Omega }h_{n}( u) \sigma_{2n} \cdot DT_{k}(
u) dx=\int_{\Omega }\sigma_{k} \cdot DT_{k} ( u) dx.
\end{align*} 
From (\ref{3.24}) and (\ref{3.25}) we get  
\begin{equation}
\limsup_{\varepsilon \rightarrow 0}\int_{\Omega }h_{n}( u^{\varepsilon
}) \aop ( x,Du^{\varepsilon }) \cdot DT_{k}(
u^\varepsilon) dx \leq \int_{\Omega }\sigma_{k} \cdot DT_{k}( u) dx.  \label{3.26}
\end{equation} 
At last, writing
\begin{multline*}
\int_{\Omega }( \aop ( x,DT_{k}( u^{\varepsilon
}) ) -\aop ( x,DT_{k}( u) ) )
\cdot
( DT_{k}( u^{\varepsilon }) -DT_{k}( u) ) dx
\\
=\int_{\Omega }\aop ( x,DT_{k}( u^{\varepsilon })
) \cdot DT_{k}( u^{\varepsilon }) dx - 
\int_{\Omega }\aop ( x,DT_{k}( u^{\varepsilon })
) \cdot DT_{k}( u) 
dx \\
-\int_{\Omega }\aop ( x,DT_{k}( u) ) \cdot (
DT_{k}( u^{\varepsilon }) -DT_{k}( u) ) dx,
\end{multline*} 
using (\ref{3.9}), (\ref{3.9bis}), (\ref{3.26}) and the monotone
character of the operator $\aop$ allow to conclude
the proof of Lemma \ref{lemme3.3}.
\end{proof}

From Lemma \ref{lemme3.3} we deduce that $\forall k>0$ 
\begin{equation*}
\lim_{\varepsilon \rightarrow 0}\int_{\Omega }\aop (
x,DT_{k}( u^{\varepsilon }) ) \cdot DT_{k}(
u^{\varepsilon }) dx=\int_{\Omega }\sigma_{k}\cdot DT_{k}(
u) dx,
\end{equation*} 
which gives thanks to a Minty argument 
\begin{equation}\label{mintydurdur}
\forall k>0,\qquad \sigma_{k}=\aop ( x,DT_{k}(
u) ) \text{ \ a.e. in }\Omega .
\end{equation} 
Using again Lemma \ref{lemme3.3} together with (\ref{3.9}) and (\ref{mintydurdur}) we conclude that 
\begin{equation*}
\forall k>0,\qquad \aop ( x,DT_{k}( u^{\varepsilon })
) \cdot DT_{k}( u^{\varepsilon }) \overset{\varepsilon
\rightarrow 0}{\rightharpoonup }\aop ( x,DT_{k}( u)
) \cdot DT_{k}( u) \text{ in }L^{1}( \Omega ) 
\text{--weak.}
\end{equation*}

{\slshape Step 4.} We now pass to the limit in the approximate problem. 
\par
Let $h$ be
an element of  $W^{1,\infty }( \mathbb{R}) $ with compact
support, let $k>0$ such that $\text{supp}( h) \subset \left[ -k,k\right] 
$ and let $\varphi $ be an element of  $W_{\Gamma _{d}}^{1,p}( \Omega
) \cap L^{\infty }( \Omega )$. Plugging the test function $ 
\varphi h( u^{\varepsilon }) $ in (\ref{3.4}) yields 
\begin{multline} \label{3.29}
\int_{\Omega }\lambda _{\varepsilon }( x,u^{\varepsilon })
\varphi h( u^{\varepsilon }) dx+\int_{\Omega }h(
u^{\varepsilon }) \aop ( x,Du^{\varepsilon }) \cdot
D\varphi dx
\\ +\int_{\Omega }\varphi h^{\prime }( u^{\varepsilon }) 
\aop ( x,Du^{\varepsilon }) \cdot Du^{\varepsilon }dx 
+\int_{\Omega }h( u^{\varepsilon }) \Phi _{\varepsilon
}( x,u^{\varepsilon }) \cdot D\varphi dx 
\\
{} +\int_{\Omega }\varphi
h^{\prime }( u^{\varepsilon }) \Phi _{\varepsilon
}( x,u^{\varepsilon }) \cdot Du^{\varepsilon }dx 
=\int_{\Omega } f ^{\varepsilon }\varphi h( u^{\varepsilon
}) dx. 
\end{multline} 
Let us pass to the limit in (\ref{3.29}) as $\varepsilon $ goes to zero.
Since $h$ has a compact support, assumption (\ref{2.4bis}) and the pointwise
convergence of $u^{\varepsilon }$ give that the field $\lambda _{\varepsilon
}( x,u^{\varepsilon }) \varphi h( u^{\varepsilon }) $
converges to $\lambda ( x,u) \varphi h( u) $ a.e. in $ 
\Omega $ and in $L^{\infty }( \Omega ) $ weak--$*$. From
(\ref{3.3}) and (\ref{3.8}) it follows that $ f ^{\varepsilon }\varphi 
h( u^{\varepsilon }) $ converges strongly to $ 
 f \varphi h( u) $ in $L^{1}( \Omega ) $. Using
assumption (\ref{2.5}) together with (\ref{3.8}) (and since $\text{supp}(
h) $ is compact) and Lebesgue's convergence theorem we obtain that $ 
h( u^{\varepsilon }) \Phi _{\varepsilon }(
x,u^{\varepsilon }) $ converges strongly to $h( u)  
\Phi( x,u) $ in $L^{N/(p-1)}( \Omega )$.
Recalling that $N\geq p$ leads to 
\begin{equation}
\lim_{\varepsilon \rightarrow 0}\int_{\Omega }h( u^{\varepsilon
}) \Phi _{\varepsilon }( x,u^{\varepsilon }) \cdot
D\varphi dx=\int_{\Omega }h( u) \Phi ( x,u)
\cdot D\varphi dx . \label{3.30}
\end{equation} 
\par
Similarly the weak convergence of $T_{k}( u^{\varepsilon }) $
yields that 
\begin{equation}\label{3.31}
\lim_{\varepsilon \rightarrow 0}\int_{\Omega }\varphi h^{\prime }(
u^{\varepsilon }) \Phi _{\varepsilon }( x,u^{\varepsilon
}) \cdot Du^{\varepsilon }dx
=\int_{\Omega }\varphi h^{\prime }( u) \Phi  ( x,u) \cdot DT_{k}( u)
dx, 
\end{equation} 
and from (\ref{3.9bis}) and (\ref{mintydurdur}) it follows that
\begin{equation}
\lim_{\varepsilon \rightarrow 0}\int_{\Omega }h( u^{\varepsilon
}) \aop ( x,DT_{k}( u^{\varepsilon }) )
\cdot D\varphi dx=\int_{\Omega }h( u) \aop (
x,DT_{k}( u) ) \cdot D\varphi dx.  \label{3.32}
\end{equation}
Finally due to Lemma \ref{lemme3.3} we have 
\begin{equation}\label{3.33}
\begin{split}
\lim_{\varepsilon \rightarrow 0}\int_{\Omega } \varphi 
h^{\prime }(u^{\varepsilon }) \aop ( x,Du^{\varepsilon })
\cdot Du^{\varepsilon }dx 
&=\lim_{\varepsilon \rightarrow 0}\int_{\Omega
}  \varphi h^{\prime }( u^{\varepsilon })\aop (
x,DT_{k}( u^{\varepsilon }) ) \cdot DT_{k}(
u^{\varepsilon }) dx
\\
&=\int_{\Omega } \varphi h^{\prime }( u) 
\aop ( x,DT_{k}( u) ) \cdot DT_{k}( u)
dx. 
\end{split}
\end{equation} 
Due to (\ref{3.29})--(\ref{3.33}) we obtain  that the field $u$ verifies
condition (\ref{2.9}) of Definition \ref{def2.2}. Condition (\ref{2.11}) is a
consequence of (\ref{2.3}), (\ref{3.9}) and Lemma (\ref{lemme3.2}).
\par
The proof
of Theorem \ref{th3.1} is now complete. 
\end{proof}

\section{Uniqueness results}

As mentioned in the introduction, we give in this section two uniqueness
results. In Theorem \ref{th4.1} below we establish that  if $\lambda ( x,s) $ is
strictly monotone with respect to $s$ then the renormalized solution of (\ref 
{1.1})--(\ref{1.3}) is unique under a  local Lipschitz condition on $\Phi
( x,s) $ with respect to $s$.

\begin{theorem}
\label{th4.1}Assume that (\ref{2.1})--(\ref{2.5BIS}) hold true. Moreover
assume that
\begin{gather}
( \lambda ( x,r) -\lambda ( x,r^{\prime })
) ( r-r^{\prime }) >0\qquad \forall r,r^{\prime }\in 
\mathbb{R},\ r\neq r^{\prime }\quad \text{a.e. }x\in \Omega;  \label{4.1}
\intertext{for any compact $C\subset \R$ there exists $L_C>0$ such that}
\left| \Phi ( x,r) -\Phi ( x,r^{\prime }) \right|
\leq L_{C}\left|  b ( x) \right| \left| r-r^{\prime
}\right| \qquad \forall r,r^{\prime }\in C.  \label{4.2}
\end{gather}
Then the renormalized solution of equation (\ref{1.1})--(\ref{1.3}) is unique.
\end{theorem}
\par

When $\lambda ( x,\cdot ) $ is assumed to be monotone
and when $1<p\leq 2$ we must replace condition (\ref{4.2}) by a global
condition and the strong monotonicity of the operator $\aop (
x,\cdot ) $ is needed.

\begin{theorem}\label{th4.2}
Assume that (\ref{2.1})--(\ref{2.5BIS}) hold true.
Moreover assume that
\begin{gather}
1<p\leq 2 < N;  \label{4.3}
\\
( \lambda ( x,r) -\lambda ( x,r^{\prime })
) ( r-r^{\prime }) \geq 0\qquad \forall r,r^{\prime }\in 
\mathbb{R},\, \quad \text{a.e. }x\in \Omega;  \label{4.3bis}
\\
( \aop ( x,\xi ) -\aop ( x,\xi ^{\prime
}) ) \cdot ( \xi -\xi ^{\prime }) \geq \alpha \frac{\left|
\xi -\xi ^{\prime }\right| ^{2}}{( \left| \xi \right| +\left| \xi
^{\prime }\right| ) ^{2-p}}\qquad \forall \xi ,\xi ^{\prime }\in 
\mathbb{R}^{N},\quad \text{a.e. }x\in \Omega;  \label{4.4}
\\
\text{there exist }L>0\text{ and }\gamma <p-\frac{3}{2}\text{ such that}  \label{4.5}
\\
\left| \Phi ( x,r) -\Phi ( x,r^{\prime }) \right|
\leq L\left| r-r^{\prime }\right| \left|  b ( x) \right|
( \left| r\right| +\left| r^{\prime }\right| +1) ^{\gamma }\quad
\forall r,r^{\prime }\in \mathbb{R},\quad \text{a.e. }x\in \Omega.  \notag
\end{gather} 
Then the renormalized solution of equation (\ref{1.1})--(\ref{1.3}) is unique.
\end{theorem}

\begin{remark} An example of function $\Phi$ verifying growth
  condition (\ref{2.5})--(\ref{2.5bis}) and (\ref{4.5}) is
  $b(x)(1+|r|)^{p-1}$ with $b$ satisfying regularity assumption
  (\ref{2.5bis}). Roughly speaking, condition (\ref{4.5}) 
  implies  that $\big|\frac{\partial\Phi(x,r)}{\partial r}\big|\leq
  |b(x)|(1+|r|)^{\gamma}$.
  When $p>3/2$ a global Lipschitz condition on $\Phi(x,r)$ with
  respect to $r$ is allowed (or a strong control of the modulus of
  continuity).  If $p\leq 3/2$ it follows that $\gamma<0$
  and then $\big|\frac{\partial\Phi(x,r)}{\partial r}\big|$ goes to
  zero as $|r|$ tends to $\infty$ almost everywhere in $\Omega$.
\end{remark}

\begin{proof}[Proof of Theorem \ref{th4.1}] Let $u$ and $v$ be two renormalized
solutions of Problem (\ref{1.1})--(\ref{1.3}). Our goal is to prove
that 
$\int_{\Omega }\left| \lambda ( x,u) -\lambda (
x,v) \right| dx=0$.
\par
Let $q>0$, $\sigma >0$ and $n\geq 1$. Using $T_{\sigma
}(T_{q}( u) -T_{q}( v) )h_{n}( u) $ which
belongs to  $W_{\Gamma _{d}}^{1,p}( \Omega ) \cap
L^{\infty }( \Omega ) $ in (\ref{2.9})$_u$ gives 
\begin{multline} \label{4.6}
\int_{\Omega }\lambda ( x,u) T_{\sigma }(T_{q}( u)
-T_{q}( v) )h_{n}( u) dx+\int_{\Omega }h_{n}(
u) \aop ( x,Du) \cdot DT_{\sigma }(T_{q}(
u) -T_{q}( v) )dx 
 \\
+\int_{\Omega }h_{n}^{\prime }( u) T_{\sigma }(T_{q}(
u) -T_{q}( v) )\aop ( x,Du) \cdot Dudx 
+\int_{\Omega }h_{n}( u) \Phi ( x,u) \cdot
DT_{\sigma }(T_{q}( u) -T_{q}( v) )dx
\\
{}+\int_{\Omega}  h_{n}^{\prime}( u)
T_{\sigma }(T_{q}( u) -T_{q}( v) ) \Phi ( x,u) \cdot Dudx 
=\int_{\Omega } f T_{\sigma }(T_{q}( u) -T_{q}(
v) )h_{n}( u) dx.  
\end{multline} 
It is then easy to pass to  the limit as $q$ tends to $+\infty $ for fixed 
$\sigma >0$ and $n\geq 1$. Indeed, since 
$\text{supp}( h_{n}) \subset \left[ -2n,2n\right] $ one has 
\begin{equation*}
h_{n}( u) DT_{\sigma }(T_{q}( u) -T_{q}( v)
)=h_{n}( u) DT_{k}(u-v)\quad \text{a.e. }x\text{ \ in }\Omega
\end{equation*} 
as soon as $q>2N+\sigma$. Moreover $T_{\sigma }(T_{q}(
u) -T_{q}( v) )$  converges to $T_{\sigma }(u-v)$  a.e. in $ 
\Omega $ and in $L^{\infty }( \Omega ) $ weak--$* $ as 
 $q$ goes to $+\infty$. Using such a process in (\ref 
{1.1})--(\ref{1.3}) written for $v$  gives by  subtraction
\begin{gather} \label{4.7}
\frac{1}{\sigma }\int_{\Omega }( \lambda ( x,u)
h_{n}( u) -\lambda ( x,v) h_{n}( v) )
T_{\sigma }(u-v)dx  
\\
+\frac{1}{\sigma }\int_{\Omega }( h_{n}(
u) \aop ( x,Du) -h_{n}( v) \aop (
x,Dv) ) \cdot DT_{\sigma }(u-v)dx 
\notag
\\
+\frac{1}{\sigma }\int_{\Omega }( h_{n}( u) \Phi  
( x,u) -h_{n}( v) \Phi ( x,v)
) \cdot DT_{\sigma }(u-v)dx  
\notag
\\
+\frac{1}{\sigma }\int_{\Omega }
 T_{\sigma}(u-v) \Big(h_{n}^{\prime }( u)
\aop ( x,Du) \cdot Du-h_{n}^{\prime }( v) \aop ( x,Dv)
\cdot Dv \Big) dx  
\notag 
\\
+\frac{1}{\sigma }\int_{\Omega }T_{\sigma }(u-v)\Big(
h_{n}^{\prime }(u) \Phi ( x,u) \cdot Du-
h_{n}^{\prime }( v) \Phi  
( x,v) \cdot Dv \Big) dx  
\notag
\\
=\frac{1}{\sigma }\int_{\Omega } f T_{\sigma }(u-v)(
h_{n}( u) -h_{n}( v) ) dx. 
\notag
\end{gather}

We  now pass to the limit successively as $\sigma \rightarrow 0$ and
as $n\rightarrow +\infty $ in (\ref{4.7}).  Since 
\begin{equation}
\frac{T_{\sigma }(u-v)}{\sigma }\overset{\sigma \rightarrow 0}{ 
\longrightarrow }\text{sign}(u-v)\indi{\{u\neq v\}}
\text{ \ a.e. in }\Omega \text{ and in }L^{\infty
}( \Omega ) \text{ weak--}* ,  \label{4.8}
\end{equation} 
and because both functions $h_{n}( u)$ and $h_n(v)$ converge to 1 almost
everywhere in $\Omega$ and in $L^{\infty }( \Omega )$ weak--$*$, 
we obtain thanks to Lebesgue's convergence theorem
\begin{equation}
\lim_{n\rightarrow \infty }\lim_{\sigma \rightarrow 0}\frac{1}{\sigma } 
\int_{\Omega }( \lambda ( x,u) h_{n}( u) -\lambda
( x,v) h_{n}( v) ) T_{\sigma
}(u-v)dx=\int_{\Omega }( \lambda ( x,u) -\lambda (
x,v) ) \text{sign}(u-v)dx  \label{4.10}
\end{equation} 
and
\begin{equation}
\lim_{n\rightarrow \infty }\lim_{\sigma \rightarrow 0}\frac{1}{\sigma } 
\int_{\Omega } f T_{\sigma }(u-v)( h_{n}( u)
-h_{n}( v) ) dx=0.  \label{4.11}
\end{equation}
\par
One has for any $\sigma>0$ and any $n\geq1$,
\begin{equation*}
\left| \frac{1}{\sigma }\int_{\Omega }h_{n}^{\prime }( u)
T_{\sigma }(u-v)\aop ( x,Du) \cdot Dudx\right| \leq \frac{1}{ 
n}\int_{\left\{ \left| u\right| <2n\right\} }\aop ( x,Du)
\cdot Dudx
\end{equation*} 
and 
\begin{equation*}
\left| \frac{1}{\sigma }\int_{\Omega }h_{n}^{\prime }(u)
T_{\sigma }(u-v) \Phi ( x,u) \cdot Dudx\right| \leq \frac{1}{n} 
\int_{\left\{ \left| u\right| <2n\right\} }\left| \Phi (
x,u) \cdot Du\right| dx.
\end{equation*} 
Therefore (\ref{2.8}) and (\ref{2.11}) imply that 
\begin{equation}
\lim_{n\rightarrow \infty }\limsup_{\sigma \rightarrow 0}\left| 
\frac{1}{\sigma }\int_{\Omega }h_{n}^{\prime }( u) T_{\sigma
}(u-v)\aop ( x,Du) \cdot Dudx\right| =0,  \label{4.12}
\end{equation} 
\begin{equation}
\lim_{n\rightarrow \infty }\limsup_{\sigma \rightarrow 0}\left| 
\frac{1}{\sigma }\int_{\Omega }h_{n}^{\prime }(u)
T_{\sigma }(u-v) \Phi ( x,u) \cdot Dudx\right| =0.  \label{4.13}
\end{equation} 
By the same way we obtain (\ref{4.12}) and (\ref{4.13}) for $v$ and
then the forth and fifth terms of (\ref{4.7}) tend to zero.
\par
To study the behavior of the second term of 
(\ref{4.7}),  we split it as follows
\begin{multline}\label{4.14}
\frac{1}{\sigma }\int_{\Omega }\big( h_{n}( u) \aop  
( x,Du) -h_{n}( v) \aop ( x,Dv)
\big) \cdot DT_{\sigma }(u-v)dx  
\\
= \frac{1}{\sigma }\int_{\Omega }h_{n}( u) \big( \aop  
( x,Du) -\aop ( x,Dv) \big) \cdot DT_{\sigma
}(u-v)dx 
\\
+\frac{1}{\sigma }\int_{\Omega }\big(
h_{n}( u) -h_{n}( v) \big) \aop ( x,Dv) \cdot DT_{\sigma }(u-v)dx.
\end{multline}
Since $h_{n}$ is a Lipschitz continuous function we have for 
$0<\sigma \leq 1$ 
\begin{multline*}
\frac{1}{\sigma }\Big|\int_{\Omega } (
h_{n}( u) -h_{n}( v) ) \aop ( x,Dv) \cdot DT_{\sigma
}(u-v)dx \Big| \\
\leq \frac{1}{n}\int_{\left\{ 0<\left| u-v\right| <\sigma
\right\} }\left| \aop ( x,DT_{2n+1}( u) ) \right|
\left| DT_{2n+1}( v) -DT_{2n+1}( v) \right| dx
\end{multline*}
which gives using  Lebesgue's convergence theorem 
\begin{equation*}
\lim_{\sigma \rightarrow 0}\frac{1}{\sigma }\Big| \int_{\Omega } 
(
h_{n}( u) -h_{n}( v) ) 
\aop ( x,Dv)  \cdot DT_{\sigma }(u-v)dx\Big| =0.
\end{equation*} 
Since $h_{n}(u)$ is non negative  the monotone character of the operator $\aop $
 and  (\ref{4.14})
 lead to $\forall n\geq 1$,
\begin{equation} \label{4.15}
\limsup_{\sigma \rightarrow 0}
\frac{1}{\sigma }\int_{\Omega }( h_{n}( u) \aop  
( x,Du) -h_{n}( v) \aop ( x,Dv)
) \cdot DT_{\sigma }(u-v)dx\geq 0.
\end{equation} 
We now deal with the third term of (\ref{4.7}). We
have 
\begin{multline} \label{4.16}
\frac{1}{\sigma }\int_{\Omega }( h_{n}( u) \Phi  
( x,u) -h_{n}( v) \Phi ( x,v)
) \cdot DT_{\sigma }(u-v)dx   \\
=\frac{1}{\sigma }\int_{\left\{ 0<\left| u-v\right| <\sigma \right\}
}h_{n}( u) ( \Phi ( x,u) -\Phi 
( x,v) ) \cdot DT_{\sigma }(u-v)dx 
\\
{} +\frac{1}{\sigma }
\int_{\left\{ 0<\left| u-v\right| <\sigma \right\} }( h_{n}(
u) -h_{n}( v) ) \Phi ( x,v) \cdot
DT_{\sigma }(u-v)dx.
\end{multline}
Since  $\text{supp}( h_{n}) = \left[ -2n,2n\right] $ we
get for $0<\sigma \leq 1$
\begin{eqnarray*}
&&\frac{1}{\sigma }\left| \int_{\left\{ 0<\left| u-v\right| <\sigma \right\}
}h_{n}( u) ( \Phi ( x,u) -\Phi  
( x,v) ) \cdot DT_{\sigma }(u-v)dx\right| \\
&&\leq \frac{1}{ 
\sigma }\int_{\left\{ 0<\left| u-v\right| <\sigma ,\left| u\right|
<2n,\left| v\right| <2n+1\right\} }h_{n}( u) \left| \Phi  
( x,u) -\Phi ( x,v) \right| \left|
Du-Dv \right| dx,
\end{eqnarray*} 
which gives thanks to assumption (\ref{4.2})
\begin{eqnarray}
&&\frac{1}{\sigma }\left| \int_{\left\{ 0<\left| u-v\right| <\sigma \right\}
}h_{n}( u) ( \Phi ( x,u) -\Phi  
( x,v) ) \cdot DT_{\sigma }(u-v)dx\right| \notag \\
&&\leq L\int_{\{ 0<\left| u-v\right| <\sigma 
\} }\left|  b ( x) \right| \big(
\left| DT_{2n+1}( u) \right| +\left| DT_{2n+1}( v)
\right| \big) dx  \label{4.17}
\end{eqnarray} 
where $L$ does not depend on $\sigma$.

Using again the fact that $h_{n}$ is Lipschitz continuous together with
assumption (\ref{2.5}) we have for
$0<\sigma \leq 1$
\begin{multline*}
\left| \frac{1}{\sigma }\int_{\left\{ 0<\left| u-v\right| <\sigma \right\}
}( h_{n}( u) -h_{n}( v) ) \Phi  
( x,v) \cdot DT_{\sigma }(u-v)dx\right|  
\\
\begin{aligned}[t]
\leq &\frac{1}{n}\int_{\{ 0<\left| u-v\right| <\sigma ,\left| u\right|
<2n,\left| v\right| <2n+1\}
}\left| \Phi ( x,v) \cdot DT_{\sigma }(u-v)\right| dx 
\\
\leq & M \int_{\{ 0<\left|
u-v\right| <\sigma \} }\left| 
 b ( x) \right| ( \left| DT_{2n+1}( u)
\right| +\left| DT_{2n+1}( v) \right| ) dx,
\end{aligned}
\end{multline*}
with $M>0$ not depending on $\sigma$.
From (\ref{4.16}) and (\ref{4.17})  it follows that for $0<\sigma \leq 1$
\begin{multline*}
\left| \frac{1}{\sigma }\int_{\Omega }( h_{n}( u)  
\Phi( x,u) -h_{n}( v) \Phi (
x,v) ) \cdot DT_{\sigma }(u-v)dx\right| 
\\
\leq (
L+M ) \int_{\{ 0<\left| u-v\right| <\sigma 
\} }\left|  b ( x)
\right| \big( \left| DT_{2n+1}( u) \right| +\left|
DT_{2n+1}( v) \right| \big) dx.
\end{multline*} 
Since $\left|b ( x) \right| ( \left| DT_{2n+1}( u)
\right| +\left| DT_{2n+1}( v) \right| ) $ lies in 
$L^{1}( \Omega )$,   Lebesgue's convergence theorem implies that
\begin{equation}
\lim_{\sigma \rightarrow 0} \frac{1}{\sigma }\int_{\Omega }(
h_{n}( u) \Phi ( x,u) -h_{n}( v) 
\Phi ( x,v) ) \cdot DT_{\sigma }(u-v)dx  =0.
\label{4.19}
\end{equation} 
Gathering (\ref{4.7}), (\ref{4.10})--(\ref{4.15}) and (\ref{4.19}) yields
\begin{equation}
\int_{\Omega }( \lambda ( x,u) -\lambda ( x,v)
) \text{sign}(u-v)dx\leq 0.  \label{4.20}
\end{equation} 
The strict monotonicity of $\lambda ( x,\cdot ) $ allows to
conclude that $\int_{\Omega }\left| \lambda ( x,u) -\lambda
( x,v) \right| dx=0$ and then $u=v$ almost everywhere in $\Omega$.
\end{proof}

\begin{remark}
In the pure Dirichlet case (i.e. $\Gamma_n=\emptyset$) if $\phi\,:\,
\R\longmapsto \R^N$  is a continuous function without any growth
assumption then there exists a renormalized solution of the problem
\begin{align}
\label{rmkuni1} 
  \null & \lambda |u|^{p-2} u - \diw(\aop(x,Du)+\phi(u))=f-\diw(g) \quad\text{in
  $\Omega$},
\\
 & u=0\quad\text{on
$\partial \Omega$},  \label{rmkuni2} 
\end{align}
with $f\in L^1(\Omega)$ and $g\in L^{p'}(\Omega)$
 (see \cite{LM} and \cite{Mur93} in
the linear case when $g\equiv 0$ ; notice that $\Gamma_n=\emptyset$ is crucial for
this existence result). When $\lambda>0$ and under a local Lipschitz
hypothesis on $\phi$, the method used in the proof of Theorem
\ref{th4.1} and the property (see \cite{BR98})
\begin{equation*}
  \diw(h_n(u)\phi(u))-h'_n(u)\phi(u)\cdot Du= \diw(\Psi_n(u)) \quad
  \text{in } \mathcal{D}'(\Omega),
\end{equation*}
where $\Psi_n(r)=\int_0^r h_n(s) \phi'(s) ds$ allow to obtain that
the renormalized solution  of (\ref{rmkuni1})--(\ref{rmkuni2}) is unique.
\end{remark}

\begin{proof}[Proof of Theorem \ref{th4.2}.] 
Let $u$ and $v$ be two renormalized solutions of problem
(\ref{1.1})--(\ref{1.3}).
\par

Let $\sigma $ be a positive real number and $n\geq 1$. Using $h=h_{n}$ and $ 
\psi =h_{n}( v) T_{\sigma }( u-v) $ in
(\ref{2.9}) written in $u$
together with similar arguments already used in the proof of
Theorem~\ref{th4.1} 
yield (by subtraction with the equivalent equality 
written in $v$)
\begin{gather} \label{4.21}
\int_{\Omega }h_{n}( u) h_{n}( v) ( \lambda
( x,u) -\lambda ( x,v) ) T_{\sigma }(
u-v) dx 
\\
 +\int_{\Omega }h_{n}( u) h_{n}( v) ( 
\aop ( x,Du) -\aop ( x,Dv) ) \cdot
DT_{\sigma }(u-v)dx  
\notag
\\
+\int_{\Omega }h_{n}( v) h_{n}( u) (  
\Phi ( x,u) -\Phi ( x,v) ) \cdot
DT_{\sigma }(u-v)dx 
\notag
\\
+\int_{\Omega }h_{n}^{\prime }( u) h_{n}( v)
T_{\sigma }(u-v)( \aop ( x,Du) +\Phi (
x,u) -\aop ( x,Dv) -\Phi ( x,v)
) \cdot Dudx 
\notag
\\
+\int_{\Omega }h_{n}^{\prime }( v) h_{n}( u)
T_{\sigma }(u-v)( \aop ( x,Du) +\Phi (
x,u) -\aop ( x,Dv) -\Phi ( x,v)
) \cdot Dvdx=0. \notag 
\end{gather} 
We now pass to the limit as $n$ goes to $+\infty $ first and then as $\sigma 
$ goes to $0$. It is worth noting that the reverse is performed in the proof
of Theorem \ref{th4.1}. Indeed passing to the limit as $\sigma \rightarrow 0$
first and then $n\rightarrow \infty $ leads to uniqueness of the solution
only in the case where $\lambda ( x,\cdot ) $ is strictly
monotone (assumption  (\ref{4.1}) in Theorem \ref{th4.1}). In the case of
Theorem \ref{th4.2}, the zero order term  namely $\lambda ( x,\cdot
)$, is monotone (see assumption (\ref{4.3bis})), and the uniqueness
proof program of Theorem \ref{th4.1} yields $\int_{\Omega }\left|
\lambda ( x,u) -\lambda ( x,v) \right| dx=0$ which is
not sufficient to ensure uniqueness. It follows that the second term of (\ref 
{4.21}) leads us to the uniqueness of the field $u$ letting first $n\rightarrow
\infty $ and then $\sigma \rightarrow 0$. This explains why global
condition and strong monotonicity  are assumed.
\par

Since $h_{n}\geq 0$ the monotone character of $\lambda ( x,\cdot ) $
implies that
 $\forall n\geq 1$
\begin{equation}
\int_{\Omega }h_{n}( v) h_{n}( u) ( \lambda
( x,u) -\lambda ( x,v) ) T_{\sigma }(
u-v) dx\geq 0.  \label{4.22}
\end{equation} 
We claim that $\forall \sigma >0$,  
\begin{equation}
\lim_{n\rightarrow \infty }\int_{\Omega }h_{n}^{\prime }( u)
h_{n}( v) T_{\sigma }(u-v)( \aop ( x,Du) +
\Phi ( x,u) -\aop ( x,Dv) -\Phi 
( x,v) ) \cdot Dudx=0.  \label{4.23}
\end{equation} 
Thanks to (\ref{2.8}) of Definition \ref{def2.2} and (\ref{2.11}) of Lemma 
\ref{lem2.4} we have 
\begin{equation}
\lim_{n\rightarrow \infty }\int_{\Omega }h_{n}^{\prime }( u)
h_{n}( v) T_{\sigma }(u-v)\aop ( x,Du) \cdot
Dudx=0  \label{4.24}
\end{equation} 
and 
\begin{equation}
\lim_{n\rightarrow \infty }\int_{\Omega }h_{n}^{\prime }( u)
h_{n}( v) T_{\sigma }(u-v)\Phi ( x,u) \cdot
Du dx=0.  \label{2.25}
\end{equation} 
Assumption (\ref{2.3}) gives with H\"{o}lder's inequality
\begin{multline*}
\left| \int_{\Omega }h_{n}^{\prime }( u) h_{n}( v)
T_{\sigma }(u-v)\aop ( x,Dv) \cdot Dudx\right|  \\
\begin{aligned}[t]
\null &\leq \frac{\beta \sigma }{n}\int_{\left\{ \left| u\right| <2n,\left|
v\right| <2n\right\} }( d( x) +\left| Dv\right|
^{p-1}) \left| Du\right| dx \\
&\leq \frac{\beta \sigma }{n}\bigg( \int_{\left\{ \left| v\right|
<2n\right\} }( \left| d( x) \right| +\left|
Dv\right| ^{p-1}) ^{\frac{p}{p-1}}dx\bigg)^{\frac{p-1}{p}}
\bigg(
\int_{\left\{ \left| u\right| <2n\right\} }\left| Du\right| ^{p}\bigg)^{ 
\frac{1}{p}} \\
&\leq \frac{\beta \sigma }{n}\left[ \Big( \int_{\Omega }\left| d 
( x) \right| ^{\frac{p}{p-1}}dx\Big)^{\frac{p-1}{p}}+ \Big(
\int_{\left\{ \left| v\right| <2n\right\} }\left| Dv\right| ^{p}dx \Big)^{ 
\frac{p-1}{p}}\right] \Big( \int_{\left\{ \left| u\right| <2n\right\}
}\left| Du\right| ^{p}dx\Big) ^{\frac{1}{p}}.
\end{aligned}
\end{multline*} 
Using again (\ref{2.8}) of Definition \ref{def2.2} together with
assumption (\ref{2.1}) leads to  
\begin{equation*}
\lim_{n\rightarrow \infty }\frac{1}{n}\int_{\left\{ \left| u\right|
<n\right\} }\left| Du\right| ^{p}dx=0,
\end{equation*} 
and since $d$ lies in $L^{p^{\prime }}( \Omega ) $
we deduce that 
\begin{equation}
\lim_{n\rightarrow \infty }\left| \int_{\Omega }h_{n}^{\prime }(
u) h_{n}( v) T_{\sigma }(u-v)\aop ( x,Dv)
\cdot Dudx\right| =0.  \label{4.26}
\end{equation} 
From assumption (\ref{2.5}) it follows that
\begin{multline*}
\left| \int_{\Omega }h_{n}^{\prime }( u) h_{n}( v)
T_{\sigma }(u-v)\Phi ( x,v) \cdot Dudx\right| 
 \\
\leq \sigma \bigg( \frac{1}{n}\int_{\left\{ \left| v\right| <2n\right\}
}  \left|  b ( x) \right| ^{\frac{p}{p-1}} \big( 1+\left|
v\right|\big) ^{p} dx\bigg)^{\frac{p-1}{p}} \bigg( \frac{1}{n} 
\int_{\left\{ \left| u\right| <2n\right\} }\left| Du\right| ^{p}dx\bigg) ^{ 
\frac{1}{p}}.
\end{multline*} 
Therefore using similar arguments to the ones used in the proof of
Lemma \ref{lem2.4} we deduce that  
\begin{equation}
\lim_{n\rightarrow \infty }\left| \int_{\Omega }h_{n}^{\prime }(
u) h_{n}( v) T_{\sigma }(u-v)\Phi (
x,v) \cdot Dudx\right| =0,  \label{4.27}
\end{equation} 
and then (\ref{4.23}) is proved.

In view of  (\ref{4.4}) and (\ref{4.5}), gathering  (\ref{4.21}),
(\ref{4.22}) and (\ref{4.23}) leads to  
\begin{multline}\label{4.28}
\alpha \int_{\Omega }h_{n}( u) h_{n}( v) \frac{\left|
DT_{\sigma }(u-v)\right| ^{2}}{ \big( \left| Du\right| +\left| Dv\right|
\big) ^{2-p}}dx 
\\
\leq \omega _{\sigma }( n) +\int_{\Omega
}h_{n}( u) h_{n}( v) \left| u-v\right| \left| b 
( x) \right| \big( \left| u\right| +\left| v\right| +1\big)
^{\gamma }\left| DT_{\sigma }(u-v)\right| dx,
\end{multline} 
where $\omega _{\sigma }( n) \overset{n\rightarrow \infty }{ 
\longrightarrow }0$. Young's inequality yields 
\begin{multline} \label{4.29}
\alpha \int_{\Omega }h_{n}( u) h_{n}( v) \frac{ 
\left| DT_{\sigma }(u-v)\right| ^{2}}{\big( \left| Du\right| +\left|
Dv\right| \big) ^{2-p}}dx  
\leq 2\omega _{\sigma }( n)
\\
{} +\frac{2\sigma ^{2}}{\alpha } 
\int_{\left\{ \left| u-v\right| <\sigma \right\}\cap \{u\neq v\} }h_{n}( u)
h_{n}( v) \left|  b ( x) \right| ^{2}\big(
\left| u\right| +\left| v\right| +1\big) ^{2\gamma }\big( \left|
Du\right| +\left| Dv\right| \big) ^{2-p}dx.  
\end{multline} 
Our goal is now to prove that $\left|  b ( x) \right|
^{2}\big( \left| u\right| +\left| v\right| +1\big)^{2\gamma }\big(
\left| Du\right| +\left| Dv\right| \big) ^{2-p}\in L^{1}( \Omega
) $ and then to pass to the limit as $n\rightarrow \infty $ in (\ref{4.29}).

If $p=2$, we have (recalling that $\text{supp}( h_{n}) =\left[ -2n,2n 
\right] $) 
\begin{multline*}
\int_{\Omega
}
h_{n}(
u) h_{n}( v) \left|  b ( x) \right|
^{2}\big( \left| u\right| +\left| v\right| +1\big) ^{2\gamma }dx 
\\
\leq \bigg( \int_\Omega
\left|  b \right| ^{N}dx \bigg)^{\frac{2}{N}}\bigg( \int_{\left\{
\left| u\right| <2n,\left| v\right| <2n\right\} } \big( \left| u\right|
+\left| v\right| +1\big) ^{\frac{2N\gamma }{N-2}}dx\bigg) ^{\frac{N-2}{N} 
}.
\end{multline*} 
Since $\gamma <p-\frac{3}{2}=\frac{1}{2}$,  we get $\frac{2\gamma N}{N-2}< 
\frac{N}{N-2}$ which implies, thanks to Lemma \ref{lem2.4}  
\begin{equation*}
\int_{\Omega }\left| u\right| ^{\frac{2N\gamma }{N-2}}dx<\infty \quad
\text{and}\quad \int_{\Omega }\left| v\right| ^{\frac{2N\gamma }{N-2}}dx<\infty.
\end{equation*} 
It follows that $\left|  b ( x) \right| ^{2}\big( \left|
u\right| +\left| v\right| +1\big) ^{2\gamma }$ $\in $ $L^{1}( \Omega
) $. 
\par
If $1<p<2$, making use of H\"{o}lder's inequality we have 
\begin{multline}\label{4.30}
\int_{\Omega }h_{n}( u) h_{n}( v) \left|  b  
( x) \right| ^{2}\big( \left| u\right| +\left| v\right|
+1\big) ^{2\gamma }\big( \left| Du\right| +\left| Dv\right| \big)
^{2-p}dx  
\\
\begin{aligned}[t]
\null & \quad \leq\bigg( \int_{\Omega }h_{n}( u) h_{n}( v) \left| 
 b \right| ^{\frac{N}{p-1}}dx\bigg) ^{\frac{2( p-1) }{N}}
\\
  &\times \bigg(  \int_{\left\{ \left| u\right| <2n,\left| v\right| <2n\right\}
}h_{n}( u) h_{n}( v) \big( \left| u\right| +\left|
v\right| +1 \big)^{\frac{2N\gamma }{N-2( p-1) }} \big(
\left| Du\right| +\left| Dv\right| \big) ^{\frac{( 2-p) N}{ 
N-2( p-1) }} dx\bigg) ^{\frac{N-2( p-1) }{N}}
\end{aligned}
\end{multline}

Let $m$ be a positive real number which be fixed in the sequel. Using
H\"{o}lder's inequality we get  
\begin{multline}\label{4.31}
\int_{\left\{ \left| u\right| <2n,\left| v\right| <2n\right\} }h_{n}
(u) h_{n}( v) \big( \left| u\right| +\left| v\right|
+1 \big)^{\frac{2N\gamma }{N-2( p-1) }} \big( \left|
Du\right| +\left| Dv\right| \big) ^{\frac{( 2-p) N}{N-2(
p-1) }} dx   
\\
\leq \bigg( \int_{\left\{ \left| u\right| <2n,\left| v\right| <2n\right\}
}h_{n}( u) h_{n}( v) \big(\left| u\right| +\left|
v\right| +1 \big)^{\nu }dx  \bigg) ^{\frac{2( p-1) (
N-p) }{p( N-2( p-1) ) }}
\\
\times \bigg( \int_{\Omega
}h_{n}( u) h_{n}( v) \frac{\big( \left| Du\right|
+\left| Dv\right| \big) ^{p}}{\big( \left| u\right| +\left| v\right|
+1\big)^{1+m}}dx \bigg) ^{\frac{N( 2-p) }{p( N-2(
p-1) ) }} 
\end{multline} 
with $\nu =\frac{N( 2\gamma p+( 1+m) ( 2-p)
) }{2( p-1) ( N-p) }$. Because $\gamma <p- 
\frac{3}{2}$ we can choose $m>0$ such that 
\begin{equation*}
\frac{N( 2\gamma p+( 1+m) ( 2-p) ) }{ 
2( p-1) ( N-p) }
<\frac{N( 2p^{2}-4p+2) }{2( p-1) (
N-p) }
<\frac{N( p-1) }{N-p}
\end{equation*} 
which gives, using Lemma \ref{lem2.4}, $\big( \left| u\right| +\left|
v\right| +1\big) ^{\nu }\in $ $L^{1}( \Omega ) $. Since 
$m>0$ from (\ref{2.durdur}) we obtain 
\begin{equation*}
\frac{\left| Du\right| ^{p}}{\big( 1+\left| u\right| \big) ^{1+m}}\in
L^{1}( \Omega ),\quad \frac{\left| Dv\right| ^{p}}{\big(
1+\left| v\right| \big) ^{1+m}}\in L^{1}( \Omega ) 
\end{equation*}
and then 
\begin{equation*}
  \frac{\big( \left| Du\right|
+\left| Dv\right| \big) ^{p}}{\big( \left| u\right| +\left| v\right|
+1\big)^{1+m}}\in L^{1}( \Omega ).
\end{equation*}
\par
Since $h_{n}( u) h_{n}( v) \overset{ 
n\rightarrow \infty }{\longrightarrow }$ $1$ almost everywhere in $\Omega$,
  Fatou's lemma, (\ref{4.30}) and (\ref{4.31}) imply that
\begin{equation*}
\left|  b ( x) \right| ^{2}\big( \left| u\right| +\left|
v\right| +1\big) ^{2\gamma }\big( \left| Du\right| +\left| Dv\right|
\big) ^{2-p}\in L^{1}( \Omega ) .
\end{equation*}

We are now in a position to pass to the limit as $n\rightarrow \infty
$ in equation (\ref{4.29}). Fatou's lemma yields that  
\begin{equation*}
\alpha \int_{\Omega }\frac{\left| DT_{\sigma }(u-v)\right| ^{2}}{\big(
\left| Du\right| +\left| Dv\right| \big)^{2-p}}dx\leq \frac{2\sigma ^{2}}{ 
\alpha }
\int_{\left\{ \left| u-v\right| <\sigma \right\}\cap \{u\neq v\} }\left|  b  
( x) \right| ^{2}\big( \left| u\right| +\left| v\right|
+1\big) ^{2\gamma }\big( \left| Du\right| +\left| Dv\right| \big)^{2-p}dx.
\end{equation*} 
Dividing the above inequality by $\sigma^2$ and taking the 
 limit as $\sigma \rightarrow 0$  gives thanks to  Lebesgue's convergence theorem
\begin{equation}
\lim_{\sigma \rightarrow 0}\frac{1}{\sigma ^{2}}\int_{\Omega }\frac{\left|
DT_{\sigma }(u-v)\right| ^{2}}{\big( \left| Du\right| +\left| Dv\right|
\big)^{2-p}}dx=0.  \label{4.32}
\end{equation}
\par
Let us consider  the field  $\frac{h_{n}( u) T_{\sigma }(u-v)}{\sigma 
}$ which belongs to $W_{\Gamma _{d}}^{1,p}( \Omega ) \cap
L^{\infty }( \Omega ) $.
From Poincar\'{e}'s inequality we have 
\begin{align*}
 \int_{\Omega }\frac{h_{n}( u) \left| T_{\sigma }(u-v)\right| }{ 
\sigma }dx &
 \leq c\int_{\Omega }\left| D\Big( \frac{h_{n}( u)
T_{\sigma }(u-v)}{\sigma }\Big) \right| dx \\
&
\leq c\bigg( \int_{\Omega }\left| \frac{T_{\sigma }(u-v)}{\sigma }\right|
\big|h_{n}^{\prime }( u)\big| \left| Du\right| dx+\int_{\Omega }\frac{ 
h_{n}( u) \left| DT_{\sigma }(u-v)\right| }{\sigma }dx \bigg)  
\\
& \leq c\bigg( \frac{1}{n}\int_{\left\{ \left| u\right| <2n\right\} }\left|
Du\right| dx+\Big( \int_{\left\{ \left| u\right| <2n,\left| v\right|
<2n+\sigma \right\} }( \left| Du\right| +\left| Dv\right| )
^{2-p}dx\Big) ^{{1/2}}
\\
& \qquad\qquad\qquad\qquad{}\times\Big( \frac{1}{\sigma ^{2}}\int_{\Omega } 
\frac{\left| DT_{\sigma }(u-v)\right| ^{2}}{\big( \left| Du\right| +\left|
Dv\right| \big) ^{2-p}}dx\Big) ^{1/2}\bigg) 
\end{align*} 
which is licit since $2-p<p$ and (\ref{2.7}) imply that both fields $\nbOne_{\left\{
\left| u\right| <2n\right\} }Du$ and $\nbOne_{\left\{ \left| v\right|
<2n+\sigma \right\} }Dv$ lie in $L^{2-p}( \Omega )$.

Letting first $\sigma \rightarrow 0$ and (\ref{4.32}) give 
\begin{equation*}
\int_{\Omega }\nbOne_{\{  u \neq v
\} }h_{n}( u) dx\leq \frac{c}{n}\int_{\left\{ \left|
u\right| <2n\right\} }\left| Du\right| dx.
\end{equation*} 
Taking the limit as $n$ goes to infinity
and using (\ref{2.8}) we conclude that
\begin{equation*}
\int_{\Omega }   \nbOne_{\left\{ u \neq  v
\right\} }dx=0
\end{equation*}
and then $u=v$  almost everywhere  in $\Omega$.
\end{proof}

\bigskip

\newcommand{\noopsort}[1]{} \newcommand{\printfirst}[2]{#1}
  \newcommand{\singleletter}[1]{#1} \newcommand{\switchargs}[2]{#2#1}
  \def\cprime{$'$}

\end{document}